\documentclass[11pt, oneside]{article}
\usepackage{amssymb,amsmath,mathrsfs}
\usepackage{geometry,tikz}
\usepackage[unicode,pagebackref=true]{hyperref}

\newcommand{\h}{h}
\newcommand{\hb}{\bar h}
\newcommand{\rrc}{f}
\newcommand{\rr}{f }
\newcommand{\z}{\rho }
\newcommand{\Z}{\mathbb{Z}}
\newcommand{\R}{\mathbb{R}}

\newcommand{\Res}{\mathop{\mathrm{Res}}}
\newcommand{\nue}{ \mu  }

\newcommand{\invfour}{ \mathcal{I} }
\newcommand{\invlift}{ Q }
\newcommand{\Ispin}{\mathscr{I}}
\newcommand{\gi}{g _I }
\makeatletter
\newcommand{\thickbar}{\mathpalette\@thickbar}
\newcommand{\@thickbar}[2]{{#1\mkern1.5mu\vbox{
  \sbox\z@{$#1\mkern-1.5mu#2\mkern-1.5mu$}%
  \sbox\tw@{$#1\overline{#2}$}%
  \dimen@=\dimexpr\ht\tw@-\ht\z@-.8\p@\relax
  \hrule\@height1\p@ % adjust for the desired rule thickness
  \vskip\dimen@
  \box\z@}\mkern1.5mu}
}
\makeatother
\title{The $L^2$ metric for hyperbolic 2-monopoles}
\author{Guido Franchetti${}^1$ and Derek Harland${}^2$
\medskip\\
{\normalsize ${}^1$Department of Mathematics, University of Pisa, Italy}\\
{\normalsize guido.franchetti@unipi.it}
\medskip\\
{\normalsize ${}^2$School of Mathematics, University of Leeds, UK}\\
{\normalsize d.g.harland@leeds.ac.uk}
}

\begin{document}
\maketitle
\abstract{It has been recently shown that the notoriously divergent $L ^2 $ metric on the moduli space of hyperbolic monopoles can be made finite by the introduction of a modified gauge fixing condition \cite{franchetti:2024}. In this paper we compute this modified $L ^2 $ metric for the mass $\tfrac{1}{2}$ 2-monopoles. The resulting metric is actually a complex non-degenerate bilinear form, which restricts to a Riemannian metric on the 4-dimensional subspace of inversion symmetric 2-monopoles. Remarkably,  we have been able to compute the $L ^2 $ form explicitly in terms of elementary functions and elliptic integrals. Its asymptotic form matches what was expected on the basis of previous results in the literature. 
}

\tableofcontents

\newpage

\section{Introduction}
Since their introduction in \cite{atiyah:1987}, hyperbolic monopoles have been the object of considerable interest due to their remarkable mathematical structure and distinctive properties as compared with their Euclidean counterparts.  Like Euclidean monopoles, hyperbolic monopoles constitute an integrable system and correspond to spectral curves \cite{murray:1996} and rational maps \cite{atiyah:1987,jarvis:1997}.  In recent years their integrability has led to the discovery of many explicit solutions \cite{manton:2014}, which are strikingly simple in comparison with the very few known explicit Euclidean monopoles \cite{braden:2021}.  Unlike Euclidean monopoles, hyperbolic monopoles enjoy a holographic property \cite{braam:1990,norbury:2001} that links them to holomorphic spheres \cite{murray:2003}.

A longstanding historical gap in the theory of hyperbolic monopoles is the lack of $L^2$ metric on their moduli space.  It was recognised early on that the $L ^2$ metric for hyperbolic monopoles is divergent \cite{braam:1989}.  This is unfortunate, because the $L^2$ metrics for Euclidean monopoles have proved very influential: they are important examples of hyperk\"ahler metrics, and they have been useful in studying both classical and quantum dynamics of monopoles.  Numerous alternative metrics have been proposed on moduli spaces of hyperbolic monopoles \cite{braam:1990,franchetti:2017,sutcliffe:2022,hitchin:1996,nash:2007,bielawski:2013,bielawski:2013a}, but all of these proposals either drastically modify or entirely avoid the $L^2$ formalism.

Recently, it has been shown \cite{franchetti:2024} that the original $L ^2 $ prescription still results in a finite bilinear form on the moduli space provided that the gauge fixing condition, usually taken to be $L ^2 $ orthogonality to the gauge group orbits, is modified. This alternative gauge-fixing condition was originally discovered using supersymmetry \cite{figueroa-ofarrill:2014}, and it is gratifying that it also solves the problem of the divergent $L^2$ metric.  It leads to an (almost) pluricomplex geometry, which is a natural generalisation of the hyperk\"ahler geometry of Euclidean monopole moduli spaces.  However, there is no guarantee that the new $L^2$ bilinear forms are real, so they are not necessarily Riemannian metrics.

In this article we compute the new $L^2$ form for the first nontrivial case, namely the moduli space of monopoles with charge $2$ and mass parameter $\tfrac{1}{2}$ (metrics for charge $1$ and any half-integer value of the mass were computed in \cite{franchetti:2024}).  This is the direct analogue for hyperbolic monopoles of Atiyah--Hitchin's famous calculation \cite{atiyah:2014} of the metric on the space of charge 2 Euclidean monopoles.  Atiyah--Hitchin's calculation was indirect: they solved the general problem of constructing hyperk\"ahler metrics with spherical symmetry, and showed that there is a unique metric whose properties match those of the monopole metric.  Unfortunately, this strategy is not available for hyperbolic monopoles, because not enough is known about the pluricomplex geometry.

Instead, we compute the $L^2$ monopole metric directly from its definition.  This is a formidable task, and we are not aware of any direct calculations of $L^2$ monopole metrics in the literature, other than the trivial case of charge 1.  The first challenge in calculating a monopole metric is to construct the monopoles explicitly.  For charge 2 and mass parameter $\tfrac12$ this can be done  using the JNR ansatz and the correspondence between instantons and monopoles \cite{manton:2014}.

The second challenge is to solve the gauge fixing condition for each tangent vector to the moduli space.  We have been able to do this using a remarkable and little-known paper of Jackiw--Rebbi \cite{jackiw:1977a}.  They explicitly constructed all harmonic spinors in the background of a JNR instanton.  Using ideas introduced in our earlier paper \cite{franchetti:2024}, we have been able to convert these into monopole zero-modes that solve the gauge-fixing condition.  This gives a basis for the tangent space to the moduli space, but we still need to match this basis with the basis generated by moduli space coordinate vector fields.  We do so using the symmetry of the monopole.

The final step is to calculate the integrals on hyperbolic space that determine the $L^2$ bilinear form on the moduli space.  Using symmetry, we reduce the number of integrals required to three. One integral is elementary but the remaining two are highly non-trivial.  Remarkably, we have been able to evaluate all three integrals in terms of elementary functions and elliptic integrals.

This gives the $L^2$ bilinear form on the eight-dimensional moduli space, which is our main result.  Our calculation shows that it is complex, not real, and that it is nondegenerate.  However, the restriction to the four-dimensional moduli space of inversion-symmetric (or centred) 2-monopoles is real, as expected \cite{franchetti:2024}, and positive definite.  So the hyperbolic counterpart to the Atiyah--Hitchin metric is still a real Riemannian metric.

We also look in detail at the asymptotic geometry of the hyperbolic monopole metric.  The natural expectation is for the asymptotic geometry to be a Kaluza--Klein metric on a torus bundle over $(H^3\times H^3)/\mathbb{Z}_2$.  Our calculation confirms this picture, except that the gauge field in the Kaluza--Klein metric is not real.  Previously, the motion of well-separated monopoles has been studied in two particular situations: with the position and electric charge of one monopole fixed \cite{gibbons:2007}; and with the pair of monopoles being centred \cite{franchetti:2023aa}.  Both metrics were derived by treating the monopoles as point-like charged particles.  In both of these situations the dynamics predicted by our metric agree with the predictions of \cite{gibbons:2007,franchetti:2023aa}, in spite of the fact that our metric is not real.

We leave further discussion of the physical and mathematical properties of this complex metric for future work.  It would be very interesting to determine whether the pluricomplex structure is integrable, and to calculate the torsion of the associated paraconformal structure on the complexified moduli space.  Doing so would help understand the properties of the $L^2$ complex metrics more generally.  While the implications of a complex metric for monopole dynamics are not immediately apparent, we note that complex metrics do still have a role to play in physics \cite{kontsevich:2021,witten:2022,benetti:2025,hull:2026}.

The remainder of this article is organised as follows.  In Section \ref{setup} we set up the notation, recall the main  results from \cite{franchetti:2024} and give a description of the charge 2 inversion symmetric monopoles moduli space both in terms of JNR data and in terms of a gauge invariant description due to Hartshorne. In section \ref{findingeigenspinors} we solve the gauge-fixing condition
using Dirac eigenspinors. Section \ref{l2metric} is devoted to the computation of the $L ^2 $ bilinear form. Its asymptotic behaviour is analysed in section \ref{asymptotics}.  We present our conclusions and directions for further research in section \ref{conclusions}.  Three appendices contain details of our calculations.

\section{Setting up the problem}
\label{setup} 
\subsection{Conventions for hyperbolic space and its spinors}
We denote by $H ^3 $ hyperbolic space with sectional curvature $-1 $. We will  work with the upper half-space model which has metric
\begin{equation}
\frac{\mathrm{d} x ^2 + \mathrm{d} y ^2 + \mathrm{d} \rho ^2 }{\rho ^2 }
\end{equation} 
with $x,y \in (- \infty , \infty )$, $\rho \in (0, \infty )$.

We denote by $S$ the spinor bundle on $H ^3 $. It is a complex rank 2 vector bundle  equipped with a symplectic structure  $( \cdot , \cdot )$ and a quaternionic structure $s \mapsto \thickbar{s} $. The symplectic form induces an isomorphism between $S$ and its dual $S ^\ast $, $s \mapsto s ^\ast :=(s , \cdot )$. In a suitable local frame, for
\begin{equation}
\label{epsidef} 
\epsilon =\begin{pmatrix}
0 & 1 \\
-1 & 0
\end{pmatrix} ,
\end{equation}
one has 
\begin{equation}
\label{spinorformulas} 
\begin{split} 
(s  , t) &= (s _1 , s _2 )\epsilon 
\begin{pmatrix}
t _1 \\ t _2 
\end{pmatrix}  
= s _1 t _2 - s _2 t _1 , \\ 
\thickbar s &= \begin{pmatrix}
\bar{s} _2  \\ - \bar{s} _1 
\end{pmatrix}  = \epsilon\begin{pmatrix}
\bar{s} _1  \\
\bar{s} _2 
\end{pmatrix} , \qquad 
s ^\ast =  (s _1 , s _2 )\epsilon 
=(- s _2 , s _1 ).
\end{split} 
\end{equation} 

For $t=t _1 + i t _2$, $ t _i  \in \mathfrak{su}(2 )$ we write $\bar t := - t ^\dagger = t _1 - i t _2 $, where $t^\dagger $ is the conjugate transpose of $t$.
We extend the quaternionic structure to $\mathfrak{sl }(2, \mathbb{C}  )$ valued spinors by setting, for $s \in S $, $t \in \mathfrak{su }(2) $,
\begin{equation}
\thickbar{s \otimes t} 
= \thickbar s \otimes \bar t
= \thickbar s \otimes (t _1 - i t _2 ).
\end{equation} 

We denote Clifford multiplication by a dot. If $(\epsilon _i )$ is an orthonormal frame on $H ^3 $ we write $\gamma _i = \epsilon  _i \cdot $ for the generators of a representation of $ \operatorname{Cl} (3) $. We follow the convention $\gamma _i \gamma _j + \gamma _j  \gamma _i = - 2 \delta _{ ij } $. We  recall that the quaternionic structure $\psi  \mapsto \thickbar\psi  $ is compatible with Clifford multiplication in the sense that
\begin{equation}
\label{barandclifford} 
\thickbar{v \cdot \psi  } =\bar v \cdot \thickbar\psi .
\end{equation}

\subsection{Definitions and properties of the $L^2$ metric}
We start by recalling some of the relevant results  in \cite{franchetti:2024}. Let $(A, \Phi )$ be an $SU (2) $  magnetic monopole on $H ^3 $, that is a solution of the Bogomolny equation $ \mathrm{d} ^A \Phi =* F $ with $A$ an $SU (2) $ connection on a trivial rank 2 vector bundle $V$ over $H ^3 $, $F$ its curvature, $\Phi$ a section of $\operatorname{End} (V) $.   We require that $V$, $\Phi$, $A$  extend to $\partial H ^2 =S ^2 _\infty $ and that the monopole abelianises at infinity. The latter condition means that  $V |_{ S ^2 _\infty } =L ^+ \oplus L ^-  $ with $L ^\pm $ an eigenbundle of $\Phi |_{ S ^2 _\infty}  $ with eigenvalue $\mp ip $ and first Chern number   $\pm n \in \mathbb{Z}   $,  and that $A |_{ S ^2 _\infty }$ splits as a sum of $U (1) $ connections on $L ^\pm $. We assume $2p \in \mathbb{Z}  $ so that a charge $n$ hyperbolic monopole is equivalent to a circle-invariant instanton with instanton number $2np $ \cite{atiyah:1987}.

A  gauge transformation is an $SU (2) $-valued section of $\operatorname{Aut} (V) $ which extends to $S ^2 _\infty $ preserving the splitting $V |_{ S ^2 _\infty }=L ^+ \oplus L ^- $. Two monopoles are (unframed) gauge equivalent if they differ by a gauge transformation. The  unframed  moduli space $ M ^{u}_{ n,p }$  is the space of monopoles  with charge $n$ and mass $p$, modulo gauge transformations.   It is a smooth manifold of dimension $4n-1 $ \cite{braam:1990}.  
Following \cite{franchetti:2024}, we call a pair of monopoles  framed gauge equivalent if they differ by a gauge transformation and their connections have the same Chern-Simons number.  The framed moduli space $M ^f _{ n,p }$ is the space of monopoles with charge $n$ and mass $p$,  modulo framed gauge transformations. The framed moduli space is an $\mathbb{R}$-principal bundle over the unframed moduli space \cite{franchetti:2024}.

The tangent space $T_\xi M^f_{n,p}$ is by definition the quotient of the space of solutions of the linearised Bogomolny equation
\begin{equation}
\label{linbog} 
\mathrm{d} ^A \phi + [a, \Phi ] = * \mathrm{d} ^A a
\end{equation}
by infinitesimal framed gauge transformations.  An infinitesimal framed gauge transformation is a solution of \eqref{linbog} of the form $(a,\phi)=(d^A\chi,[\Phi,\chi])$ satisfying
\begin{equation} 
\label{framedinfinitesimal} 
- \frac{1}{2} \int _{ H ^3 } \operatorname{Tr} \left(  \mathrm{d} ^A \chi  \wedge *  \mathrm{d} ^A \Phi \right) =0.
\end{equation}
The complexified tangent space $T^\mathbb{C}_\xi M^f_{n,p}$ is defined similarly, with $a,\phi,\chi$ being $\mathfrak{sl}(2,\mathbb{C})$- (rather than $\mathfrak{su}(2)$-) valued.  It was shown in \cite{franchetti:2024} that $T^\mathbb{C}_\xi M^f_{n,p}$ is naturally isomorphic to the space of solutions of \eqref{linbog} and the modified gauge fixing condition
\begin{align}
\label{gaugefixingminus} 
* \mathrm{d} ^A * a + [ \Phi , \phi ] =-   2i \phi .
\end{align}
In other words, given any solution $(a,\phi)$ of \eqref{linbog}, there exists a unique infinitesimal framed  $\mathfrak{sl}(2,\mathbb{C})$ gauge transformation $\chi$ such that $(a+d^A\chi,\phi+[\Phi,\chi])$ solves \eqref{linbog} and \eqref{gaugefixingminus}.

Equation (\ref{gaugefixingminus}) differs from the standard choice by  the $- 2i \phi $ term at the right hand side. An important reason for this choice is the fact that it makes the $L ^2 $ metric on $M ^f _{ n,p }$
\begin{equation}
\label{L2product} 
g ((a , \phi ), (a, \phi )) = - \frac{1}{2} \int _{ H ^3 } \operatorname{Tr} (a \wedge * a + \phi * \phi )
\end{equation} 
finite.   However, since the gauge fixing condition is complex, there is no guarantee that the resulting bilinear form on $T_\xi M^f_{n,p}$ is real, let alone positive definite. In fact,  the gauge fixing condition with the opposite sign,
\begin{equation}
\label{gaugefixingplus} 
* \mathrm{d} ^A * a + [ \Phi , \phi ] =   2i \phi  ,
\end{equation} 
works just as well   and also leads to a finite bilinear form which is the complex conjugate of the previous one  \cite{franchetti:2024}.  Since the treatment for either sign choice is completely parallel, in this paper we focus on the choice (\ref{gaugefixingminus}).   We will say that $(a, \phi )$ is real if it is $  \mathfrak{su }(2) $-valued (that is, anti-hermitian) up to an infinitesimal gauge transformation, i.e.
\begin{equation}
a + a ^\dagger =i \mathrm{d} ^A \Lambda , \quad \phi + \phi ^\dagger =i [ \Phi , \Lambda ]
\end{equation} 
with $\Lambda : H ^3 \rightarrow \mathfrak{su }(2) $. Similarly we say that $(a, \phi )$ is purely imaginary if it is $i \, \mathfrak{su}(2) $-valued up to an infinitesimal gauge transformation, and complex otherwise.
It was shown in \cite{franchetti:2024} that $g$ is real on the unframed moduli space of charge 1 monopoles, where it becomes the product metric on $H ^3 \times \mathbb{R}  $, and on the moduli space of inversion symmetric monopoles.

Furthermore \cite{franchetti:2024}, for $\xi  =(A , \Phi ) \in M ^f _{ n,p }$  there is an isomorphism 
\begin{equation}
\label{isotangentspinors} 
T_\xi  ^ \mathbb{C}M ^f_{n,p} \simeq K^-  \otimes E ^-_\xi ,
\end{equation}
 where $K ^- $ is the space of  Killing spinors $\psi$ on $H ^3 $ with constant $- \tfrac{i}{2} $, which for any $X \in \mathfrak{X}  (M) $ satisfy the equation
\begin{equation}
\label{killingspinoreqminus} 
\nabla _X \psi = - \frac{i}{2} X \cdot \psi ,
\end{equation} 
and $E _\xi  ^- $ is the vector bundle consisting  of $\mathfrak{sl} (2, \mathbb{C}  )$ valued eigenspinors in the adjoint representation with eigenvalue $- \tfrac{i}{2}$, 
\begin{equation}
\label{eigenspinoreqminus} 
D ^A \nu - [ \Phi , \nu ]  = - \frac{i}{2} \nu .
\end{equation} 
The Killing spinor equation on $H ^3 $ has two linearly independent solutions, so $K ^\pm \simeq \mathbb{C}  ^2 $. The vector bundles $E _\xi ^\pm $ have complex rank  $2n $ \cite{franchetti:2024}.

The isomorphism (\ref{isotangentspinors})  makes use of the identification $ (\Lambda ^0 \oplus \Lambda ^1)\otimes \mathbb{C}  \simeq \operatorname{End} (S ) $, where $S$ is the spinor bundle on $H ^3 $, induced by Clifford multiplication. Let $\{ \psi _1 , \psi _2\}  $ be a basis of $K ^- $,  $(a, \phi ) $ be a solution of (\ref{linbog}), (\ref{gaugefixingminus}), $ \{ a _k \} $ be the components of $a$ with respect to some orthonormal frame $ \{ \epsilon  _k \} $ on $H ^3 $. If we define $\nu _1 $, $\nu _2 $ via
\begin{equation}
\label{tvectospinors} 
 a _i \gamma _i + \phi \operatorname{Id} = \nu _1 \psi _1 ^\ast + \nu _2 \psi _2 ^\ast
\end{equation} 
then they satisfy (\ref{eigenspinoreqminus}).  Viceversa if $\nu \in E_\xi ^- $, $\psi \in K ^- $ and we define
\begin{equation}
\label{spinorstotv} 
\phi =\frac{1}{2} \operatorname{Tr} (\nu \otimes \psi ^\ast ), \quad 
a _k  =- \frac{1}{2} \operatorname{Tr} (\nu \otimes \psi ^\ast \gamma _k )
\end{equation}  
then $(a, \phi )$ satisfy (\ref{linbog}), (\ref{gaugefixingminus}).

Denoting by $K ^+ $ the space of Killing spinors with constant $ + \tfrac{i}{2}$,
\begin{equation}
\label{killingspinoreqplus} 
\nabla _X \psi =  \frac{i}{2} X \cdot \psi ,
\end{equation} 
and by  $E ^+ _\xi  $ the space of eigenspinors with eigenvalue $+ \tfrac{i}{2} $,
\begin{equation}
\label{eigenspinoreqplus} 
D ^A \nu - [ \Phi , \nu ]  =  \frac{i}{2} \nu ,
\end{equation} 
there is a completely analogous decomposition $T_\xi^{ \mathbb{C}  } M ^f _{ n,p }\simeq K^+  \otimes E_\xi ^+ $ corresponding to the  choice (\ref{gaugefixingplus})  for the gauge fixing condition.

The quaternionic structure $\nu \mapsto \thickbar \nu $ on $S$ maps $ K ^\pm \rightarrow K ^\mp $ and $E_\xi ^\pm \rightarrow E_\xi ^\mp $. Because of the isomorphism (\ref{isotangentspinors}) if  
$(a , \phi )$ is a tangent vector then
\begin{equation}
\label{baraphi} 
( a, \phi )\mapsto (\bar{a} , \bar \phi )\quad \Leftrightarrow \quad \nu \otimes \psi ^\ast \mapsto \thickbar \nu \otimes \thickbar \psi ^\ast .
\end{equation} 
For  $u=(a , \phi )$  a tangent vector we will use the notation  
\begin{equation}
\label{barontangentspace} 
\bar u :=( \bar{a} , \bar \phi ).
\end{equation}

\subsection{The moduli space of 2-monopoles}
\label{invsym2monopoles} 
We now specialise to the case of inversion symmetric monopoles with charge $n =2 $ and mass $p =\tfrac{1}{2}$. Any such monopole is equivalent to an anti self-dual (ASD) circle-invariant instanton $B $ on Euclidean 4-space $E ^4 $ \cite{atiyah:1987} coming from the JNR  ansatz  \cite{jackiw:1977a}
\begin{equation}
\label{cijnransatz} 
B =B _a  \, \mathrm{d} x ^a , \quad B _a = \frac{i}{2} \bar \sigma_{ a b }\partial ^b \log f.
\end{equation} 
Here  $x ^a =(x,y,z,w )$ are Cartesian coordinates on $E ^4 $, $\bar\sigma _{ a b  }= - \bar \sigma _{ b a }$ is self-dual, $\bar\sigma _{ ij } = \epsilon _{ ij k } \sigma _k $, $\bar\sigma _{ i4 } = \sigma _i $, with $\{ \sigma _1 , \sigma _2 , \sigma _3 \} $ the Pauli matrices, and $f$ is a circle-invariant harmonic function on $E ^4 $ with three distinct poles $\{ p _i \} $.   Let us take the circle to act on the $(z,w )$ plane and introduce polar coordinates $(\rho , \theta )$ on that plane.  In order for $f$ to be circle invariant its poles have to lie in the $(x,y )$-plane, so
\begin{equation}
\label{rrc} 
f(x) = \sum _{i =1 }^{3}  \frac{\lambda _i ^2  }{(x- x _i )^2 + (y- y _i ) ^2 + \rho  ^2  }.
\end{equation}
Without loss of generality we can rescale $f$ by a constant so from now on we take $\lambda _3 =1 $. The resulting connection $B$ in \eqref{cijnransatz} will be circle-invariant up to a gauge transformation.  It is possible to choose a gauge in which $B$ is strictly circle-invariant.  In this gauge we may write $B =A - \Phi \, \mathrm{d} \theta$, with $A (\partial _\theta )=0 $. The pair $(A , \Phi )$ is then the monopole associated to $B$.

Apparently the JNR ansatz has 8 degrees of freedom:  the positions $p _i =(x _i , y _i ,0,0)$ of the poles in the plane fixed by the circle action, and the scales $\lambda _1 $, $\lambda _2 $  of two  poles. However if the poles all lie on the same line or circle, and in no other case, the JNR ansatz produces gauge equivalent instantons. It is then possible to shift the poles along the circle or line (the shift is weight-dependent) via a gauge transformation \cite{jackiw:1977a}. Since three points always lie on a line or circle, for $n =2 $  the JNR ansatz parametrises a 7-dimensional space which is the unframed moduli space  $M ^u_{2,1/2} $. In the following we will write $M ^u $ for $M ^u _{ 2 , 1/2 }$ and $M ^f $ for $M ^f _{ 2 , 1/2 }$.

The isometry group $SO(3,1)$ of $H^3$ acts naturally on the moduli space $M$ of hyperbolic monopoles.  It is known that every charge 2 monopole can be transformed using isometries to a monopole that is invariant up to gauge transformations under the action of the dihedral group $D_2\subset SO(3,1)$ \cite{norbury:2006}.  So without loss of generality we now assume that our monopole is invariant under $D_2$.

To understand what JNR data corresponds to a $D_2$-symmetric monopole it is convenient to introduce a description of $M ^u $ in terms of  gauge invariant data due to Hartshorne \cite{hartshorne:1978}. The geometry is clearer if we work on $S ^4 $ rather than $E ^4 $. This is possible since, due to the conformal invariance of the ASD equations, an instanton on $E ^4 $ is also an instanton on $S ^4 $, and the two pictures are related by stereographic projection. 
 \begin{figure}[htbp]
\begin{center}
 \begin{tikzpicture}[scale=1.7]
\filldraw[black] (0.58,0.82) circle (2pt) node[anchor=south]{ $p _1 $};
\filldraw[black] (-0.58,0.82) circle (2pt)node[anchor=south]{ $p _2 $};
\filldraw[black] (0,-1) circle (2pt) node[anchor=north]{ $p _3 $};
\filldraw[black] (0,0.82) circle (1pt)node[anchor=north,scale=0.8]{ $a _3 $};
\filldraw[black] (0.12,-0.65) circle (1pt)node[anchor=west,scale=0.8]{ $a _2 $};
\filldraw[black] (-0.12,-0.65) circle (1pt)node[anchor=east,scale=0.8]{ $a _1 $};
\filldraw[gray] (0.,0.) circle (1pt);
\draw (0,-1) -- (0.58,0.82) ;
\draw (0,-1) -- (-0.58,0.82) ;
\draw (0.58,0.82) -- (-0.58,0.82) ;
\draw [red](0,0) -- (0,-0.82) ;
\draw [blue](0,0) -- (0.18,0) ;
\node[blue, scale=0.6] at (0.09,0.08) {1-b};
\node[red, scale=0.6] at (-0.06,-0.4) {b};
\draw (0,0) circle (1);
\draw(0,0) ellipse (0.18 and 0.81);
\end{tikzpicture}
\caption{Gauge invariant data for an inversion symmetric 2-instanton.}
\label{ponceletellipse} 
\end{center}
\end{figure}
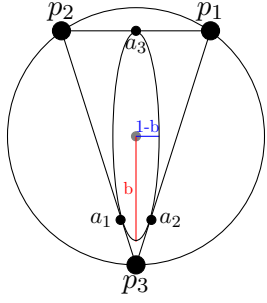
The gauge invariant data consists of a pair of conics: a circle, obtained as the intersection between $ S ^4 $ and a plane $\Pi $, and an ellipse in $\Pi $ contained in the circle. The circle and ellipse need to satisfy the conformally invariant condition  that there is  a triangle  inscribed in the circle and circumscribing  the ellipse. If a given pair of conics admits one such triangle then it admits a Poncelet porism, that is a 1-parameter family of them, and one of the triangle vertices can be placed at any point of the circle.

Let us see when a configuration is $D_2$-symmetric. Take Cartesian coordinates $\{ x _i  \} $  on $E ^5 $ and embed  $S ^4 \subset E ^5 $ in the standard way. The $E ^4 $ coordinates are obtained via stereographic projection,
\begin{equation}
\label{S4coordinates} 
(x,y,z,w ) =\frac{1}{1- x _5  } (x _1 , x _2 , x _3 , x _4 ).
\end{equation}
The circle action on $S^4$ is given by rotations of the $x_3,x_4$ coordinates (corresponding to rotations of the coordinates $z,w$ on $E^4$).  The action of $SO(3)$ on the coordinates $x_1,x_2,x_5$ gives an isometric action of $SO(3)$ on $H^3$.  The action of $D_2\subset SO(3)$ on $S^4$ is given by
\begin{equation}\label{D2onS4}
\begin{split} 
R _1 (x _1 , x _2 ,x_3,x_4, x _5 )&= (x _1 , - x _2 ,x_3,x_4, - x _5 ),\\
R _2 (x _1 , x _2 ,x_3,x_4, x _5 )&= (-x _1 ,  x _2 ,x_3,x_4, - x _5 ),\\
R _3 (x _1 , x _2 ,x_3,x_4, x _5 )&= (-x _1 , - x _2 ,x_3,x_4,  x _5 ).
\end{split} 
\end{equation}
The Poncelet circle and ellipse must be invariant under the action of the circle and of $D_2$.  Circle-invariance means that the circle and ellipse are contained in the subspace $\{x_3=x_4=0\}\subset E^5$.  Then dihedral invariance means in addition that they are contained one of the following planes: $\{x_1=x_3=x_4=0\}$, $\{x_2=x_3=x_4=0\}$, or $\{x_3=x_4=x_5=0\}$.  These planes are related by $SO(3)$ rotations, so without loss of generality we assume that they are contained in the plane $\{x_2=x_3=x_4=0\}$.

Figure \ref{ponceletellipse} shows the unit circle and a $D_2$-invariant ellipse satisfying the Poncelet condition.  In this picture, $x_1$ is the horizontal coordinate and $x_5$ is the vertical coordinate.  The Poncelet condition implies that the lengths of the semi-major and semi-minor axes sum to 1, so we denote these by $b\in(0,1)$ and $1-b$.  A $\pi/2$ rotation sends $b$ to $1-b$, so without loss of generality we take $b\in[\frac12,1)$.  When $b=\frac12$ the ellipse becomes a circle and the hyperbolic monopole has enhanced $O(2)$ symmetry.

Figure \ref{ponceletellipse} also shows a triangle satisfying the Poncelet condition.  The coordinates $(x_1,x_5)$ of the vertices $p_1,p_2,p_3$ are 
\begin{equation}
\label{refconfigcircle} 
p _3 =(0,-1), \quad 
p _1 =( \sqrt{ 1- b ^2 } ,b), \quad 
p _2 =( -\sqrt{ 1- b ^2 } ,b) .
\end{equation} 
This set of points is fixed by the rotation $R_3$, but under $R_2$ it gets mapped to
\begin{equation}
p _3' =(0,1), \quad 
p _1' =( -\sqrt{ 1- b ^2 } ,-b), \quad 
p _2' =( \sqrt{ 1- b ^2 } ,-b),
\end{equation} 
which gives another triangle in the porism.

Under stereographic projection, the points $p_1,p_2,p_3$ map to the points $(x,y,z,w)=(x_i,0,0,0)$ with
\begin{equation}
\label{polexaxis} 
x _1  =x _0:=\sqrt{\frac{1+b}{1-b}} , \quad x _2 =- x _0 , \quad x _3 =0.
\end{equation}
The constraint $b\in[\frac12,1)$ is equivalent to $x_0\in[\sqrt{3},\infty)$.  These points will be the poles of our JNR function \eqref{rrc}.  The value $x _0 =\sqrt{ 3 }$, corresponding to $b =\tfrac{1}{2}$, is the one for which the monopole has enhanced $O (2) $ symmetry.

We determine the weights $\lambda_i$ using a method of \cite{atiyah:1993}.  Let $c(x)$ be the degree $3$ polynomial with roots $x_1,x_2,x_3$:
\begin{equation}
c(x)=(x-x_1)(x-x_2)(x-x_3)=x(x^2-x_0^2).
\end{equation}
Similarly, let $(x_i',0,0,0)$ be the images of $p_i'$ under stereographic projection,
\begin{equation}
x _1'  =-x _0^{-1}, \quad x _2' = x _0^{-1} , \quad x _3' =\infty,
\end{equation}
and let $c'(x)$ be the corresponding polynomial
\begin{equation}
c'(x)=(x^2-x_0^{-2}).
\end{equation}
(This polynomial is quadratic rather than cubic because the root $x_3'$ lies at $\infty$).  The roots $x_i(t)$ of the polynomial
\begin{equation}
c(x)+tc'(x)=x(x^2-x_0^2)+t(x^2-x_0^{-2})
\end{equation}
correspond to vertices $p_i(t)=(2x_i(t),x_i(t)^2-1)/(x_i(t)^2+1)$ of triangles in the Poncelet porism.  The JNR weights satisfy
\begin{equation}\label{scales2}
\frac{\lambda_i^2}{\lambda_j^2}=\lim_{t\to 0}\frac{|p_i(0)-p_i(t)|}{|p_j(0)-p_j(t)|}.
\end{equation}
Solving the cubic perturbatively gives
\begin{equation}
x _1(t) =x _0 + \frac{t}{2} \left( \frac{1}{x _0 ^4 } -1 \right) , \quad 
x _2(t) =-x_0 + \frac{t}{2} \left( \frac{1}{x _0 ^4 } -1 \right) , \quad 
x _3(t) = - \frac{t}{ x _0 ^4 }
\end{equation}
up to terms quadratic in $t$.
Using (\ref{scales2}) and choosing $\lambda_3=1$ gives the instanton scales
\begin{equation}
\label{lambdapolesline} 
\lambda ^2 := \lambda _1 ^2 =\lambda _2 ^2 = \frac{x _0 ^4 -1 }{2}.
\end{equation} 

In later computations we will take a reference configuration of this type. More precisely we will consider the JNR data
\begin{equation} 
\label{refconfigjnr} 
\begin{split}
(x_1,y_1) &=(x _0 , 0 ), \quad (x_2,y_2) =(- x _0 ,0 ), \quad (x_3,y_3) =(0,0), \\
\lambda _1 &= \lambda _2 =\sqrt{\frac{x _0 ^4 -1 }{2}}, \quad \lambda _3 =1 . 
\end{split} 
\end{equation}
The framed moduli space $M^f$ consists of the orbit of this 1-parameter family of monopoles under the action of $SO(3,1)\times\R$, in which $SO(3,1)$ acts on $H^3$ and $t\in\R$ acts by gauge transformation with $\exp(t\Phi)$.  It is a cohomogeneity-one manifold of the form
\begin{equation}\label{framed moduli space}
M^f\cong M^u\times\R,\quad
M^u\cong \{\sqrt{3}\}\times SO(3,1)/O(2) \cup (\sqrt{3},\infty)\times SO(3,1)/D_2,
\end{equation}
in which $M^u$ is the unframed moduli space.

The $D_2$-symmetric monopoles are also invariant under an inversion.  Consider the orientation-preserving isometry
\begin{equation}
\label{invfoursphere} 
(x_1 , x _2 , x _3 , x _4 , x _5 ) \mapsto  ( - x _1 , - x _2 ,  x _3 ,- x _4 , -x _5 )
\end{equation}
of $S^4$.  It is clear that our chosen ellipse in the plane $x_2=x_3=x_4=0$ is invariant under this transformation, so the instanton is also invariant.  The corresponding transformation of $E^4$ is
\begin{equation} 
\label{invcondE4} 
(x,y,z,w ) \mapsto  \frac{1}{x ^2 + y ^2 + z ^2 + w ^2  } (-x,-y,z, -w)
\end{equation}
and the corresponding isometry of $H^3\times S^1$ is
\begin{equation}
\label{fourdiminvdef} 
 \invfour (x,y,\rho , \theta ) = (I (x,y,\rho),   - \theta),
 \end{equation}
where
\begin{equation}
\label{inv3uhs} 
I (x,y, \rho ) = \frac{1}{x ^2 + y ^2 + \rho ^2 }  (-x,-y, \rho )
\end{equation}
is the inversion of $H^3$ about the point $(x,y,\rho)=(0,0,1)$.

Since $\invfour $ is orientation-preserving, it maps ASD instantons into ASD instantons. Because of the reversal of orientation on the $S ^1 $ factor, in terms of the monopole fields we have $\invfour ^\ast (A - \Phi \mathrm{d} \theta ) =(I ^\ast A +  I ^\ast \Phi  \, \mathrm{d} \theta )$, that is  $(A , \Phi )\mapsto (I ^\ast A , - I ^\ast \Phi )$.   So a $( A , \Phi )$ is inversion symmetric if there is an $SU (2) $ gauge transformation $\gi $ such that  
\begin{equation}
\label{invsymmonopole} 
I ^\ast A = \operatorname{Ad} _{ \gi } A + \gi \mathrm{d} \gi ^{-1} ,  \quad  I ^\ast (\Phi ) =  -\operatorname{Ad} _{ \gi } \Phi .
\end{equation}
It can be proved that the moduli space $M^{I}\subset M^u$ of inversion-symmetric monopoles is the orbit of the family \eqref{refconfigjnr} of $D_2$-invariant monopoles under the action of $SO(3)\subset SO(3,1)$.   Hence it has the form
\begin{equation}
\label{modulispacecentred} 
M ^I \cong\{\sqrt{3}\}\times SO(3)/O(2) \cup (\sqrt{3},\infty)\times SO(3)/D_2.
\end{equation}

A tangent vector $(a, \phi )$ to $M^u$ is tangent to $M^I$ if, up to gauge transformations, $ I ^\ast a =a $ and $I ^\ast \phi =- \phi $.

\section{Explicit Dirac eigenspinors}
\label{findingeigenspinors} 
The aim of this Section is to find a basis for the space of solutions of \eqref{linbog}, \eqref{gaugefixingminus}.  In other words, we will find a basis for the eight-dimensional tangent space $T^{ \mathbb{C}  } _\xi  M  ^f $ that satisfies our gauge-fixing condition.  We do so by finding four linearly independent solutions of the eigenspinor equation \eqref{eigenspinoreqminus}, which form a basis for $E_\xi^-$, and applying the isomorphism \eqref{isotangentspinors}.

In fact, two solutions of \eqref{eigenspinoreqminus} are already at hand for, as  noted in \cite{franchetti:2024}, if $\psi \in K ^- $ then $F   \cdot \psi \in E_\xi^- $.  We obtain a third solution by exploiting a relationship \cite{franchetti:2024} between solutions of \eqref{eigenspinoreqminus} and harmonic spinors on $E^4$, reviewed in Section \ref{soleigenspinoreq}, together with an explicit construction \cite{jackiw:1977a} of harmonic spinors on $E^4$.  We obtain a fourth making use of the anti-involution $\nu \mapsto \Ispin (\thickbar \nu )$ of $E _\xi  ^- $ obtained by combining the map $\Ispin: E ^\mp  _\xi  \rightarrow E ^\pm  _\xi  $ induced by inversion, defined in Section \ref{sn:inversionaction}, with the action $\nu \mapsto \thickbar \nu $ of the quaternionic structure.

\subsection{Eigenspinors  from harmonic spinors on $E^4$}
\label{soleigenspinoreq} 
Let us first explain how solutions of \eqref{eigenspinoreqminus} are related to harmonic spinors on $E^4$.  As it is well known,  $E ^4 \setminus E ^2 $ is conformally equivalent to $H ^3 \times S ^1 $. Let us write $g _E $ for the metric on $E ^4 $, $g _H $ for the product metric on $H ^3 \times S ^1 $ and $\Lambda$ for the conformal factor so that
\begin{equation}
\label{confequiv} 
g _E =\Lambda ^2 g _H .
\end{equation} 
In particular if we take $S ^1 $ to act in the $(z,w) $-plane and introduce polar coordinates $(\rho , \theta )$ in that plane then $\Lambda=\rho$ and (\ref{confequiv}) takes the form
\begin{equation}
\label{hmetricuhs} 
\mathrm{d} x ^2 + \mathrm{d} y ^2 + \mathrm{d} z ^2 + \mathrm{d} w ^2 
= \rho ^2  \left(  \frac{\mathrm{d} x ^2 + \mathrm{d} y ^2 + \mathrm{d} \rho ^2 }{\rho ^2 } + \mathrm{d} \theta ^2 \right) .
\end{equation} 

Let  $\partial _\theta $ be the generator of $S ^1 $, $ (A , \Phi )$  a monopole on $H ^3 $, $B$  the pullback of $(A , \Phi )$ to $H ^3 \times S ^1  $,  $B  =A - \Phi \,   \mathrm{d} \theta $. The 1-form $B$ defines a circle-invariant connection on $H ^3 \times S ^1 $ which has, as a consequence of the Bogomolny equation, anti self-dual curvature, $* _H F ^B =- F ^B $  with $F ^B =\mathrm{d} B + B \wedge B $ and $* _H $  the Hodge star with respect to $g _H $. In other words $B$ is a circle-invariant anti self-dual instanton on $H ^3 \times S ^1 $. Since $( A , \Phi )$ extends to $S ^2 _\infty $ and  the Hodge operator is conformally invariant in middle dimension, $B$ is also a circle-invariant anti self-dual instanton on $E ^4 $.   Let  $\nue $ be a 4-spinor on $E ^4 $. By the conformal invariance of the kernel of the Dirac operator \cite{lawson:1989}
\begin{equation}
\label{harmtoeigen} 
D ^B _E \nue  =0 \quad \Leftrightarrow \quad D ^B _H (\Lambda ^{ 3/2 } \nue )=0,
\end{equation} 
where $D ^B _E $, $D ^B _H $ are the Dirac operators twisted by $B$ on $E ^4 $ and $H ^3 \times S ^1 $.  In local orthonormal frames $\{ e _k \} $ for $TE^4$ and $\{ \epsilon _k=\Lambda^{-1}e_k \} $ for $TH ^3 \times S ^1 $, they take the form
\begin{equation}
\begin{split} 
D ^B _E &= \Gamma ^k \big((\nabla _{E}) _{  e_k}   + [B (e _k ), \cdot ] \big),\quad 
D ^B _H = \Gamma ^k \big((\nabla _{H}) _{  \epsilon _k}   + [B (\epsilon  _k ), \cdot ] \big),
\end{split} 
\end{equation} 
in which  $\{ \Gamma ^1, \Gamma ^2 , \Gamma ^3  , \Gamma ^4  \} $ are generators of the Clifford algebra and $\nabla _E $, $\nabla _H $ the lifts of the Levi-Civita connections to the spinor bundle. The Clifford volume element $\Gamma  _1 \Gamma  _2 \Gamma  _3 \Gamma  _4 $ squares to 1. We say that a spinor  has positive (respectively negative) chirality if it belongs to the $+1 $ ($-1$) eigenspace. We will restrict attention to harmonic spinors $\nue  $ with positive chirality (by the Lichnerowicz formula, there are no normalisable harmonic spinors with negative chirality).

The space of harmonic spinors decomposes into weight spaces for the circle action.  Let us assume that $\tilde{\nu}=\Lambda^{3/2}\nue$ has weight $\pm\tfrac{i}{2}$.  Then
\begin{equation}
\label{euclspinorweightp} 
\nu :=e^{\mp\frac{i}{2}\theta}\tilde{\nu}=\Lambda ^{ 3/2 } e^{\mp\frac{i}{2}\theta}\nue.
\end{equation}
has weight $0$, so descends to a spinor on $H^3$.  We choose a frame for $T(H^3\times S^1)$ such that $\epsilon _4 = \partial _\theta $ and $\epsilon_1,\epsilon_2,\epsilon_3$ are a frame for $TH^3$, pulled back to $H^3\times S^1$.  Then
\begin{equation}
D^B_H\nu = \big[D^B_H,e^{\mp\frac{i}{2}\theta}\big]\tilde{\nu} + e^{\mp\frac{i}{2}\theta}D^B_H\tilde\nu=\mp\tfrac{i}{2}\Gamma_4\nu
\end{equation}
since $D^B_H\tilde{\nu}=0$.  We wish to reexpress this in three-dimensional language. There is an isomorphism $ \operatorname{Cl} ^0 (4) \rightarrow \operatorname{Cl} (3) $,  $ \Gamma ^i \mapsto \Gamma ^i \Gamma ^4 =: \gamma ^i , i =1,2,3$.  So multiplying with $-\Gamma_4$ gives
\begin{equation}\label{eigenspinoreqderiv} 
\mp\tfrac{i}{2}\nu=-\Gamma_4D^B_H\nu = \sum_{i=1}^3\gamma_i((\nabla_H)_{\epsilon_i}\nu+[A(\epsilon_i),\nu])+\Gamma_4^2[\Phi,\nu]=D^A_H\nu-[\Phi,\nu].
\end{equation}
Thus solutions of the eigenspinor equations (\ref{eigenspinoreqminus}), (\ref{eigenspinoreqplus}) correspond to harmonic spinors on $E^4$ with weight $\pm\tfrac{i}{2}$.

\subsection{Solving the eigenspinor equation}
\label{solvingtheeigenspinoreq} 
At this point it is convenient to fix some of our choices. We take $(A , \Phi )$ to be the monopole described by the JNR data corresponding to the reference configuration (\ref{refconfigjnr}), with JNR function 
\begin{equation}
\label{jnrrefconfig} 
f (x,y, \rho ) =  \lambda ^2 \left( \frac{ 1 }{(x- x _0 )^2 + y ^2 + \rho ^2 } +  \frac{1 }{(x +  x _0 )^2 + y ^2 + \rho ^2 }\right)  + \frac{1}{x ^2 + y ^2 + \rho ^2 },
\end{equation} 
where $\lambda$ is  given in terms of $x _0 $ by (\ref{refconfigjnr}). The corresponding instanton $B$ is (\ref{cijnransatz}) and, after switching to polar coordinates $(\rho , \theta )$ in the $(z,w )$-plane and to a circle-invariant gauge, the corresponding monopole $(A , \Phi )$ is defined by $B =A - \Phi \, \mathrm{d} \theta $, $A (\partial \theta  )=0 $. The computation is outlined in Appendix \ref{jnrmonopoles}. 

We take the chiral representation of  $\operatorname{Cl} (4) $ given by
\begin{equation}
\Gamma _k =\begin{pmatrix}
0 & \sigma _k  \\
- \sigma _k & 0
\end{pmatrix} ,\, k =1,2,3,  \qquad 
\Gamma _4 =\begin{pmatrix}
0 & i \\
i & 0
\end{pmatrix} ,
\end{equation} 
so that
\begin{equation}
\Gamma _k \Gamma _4 =\begin{pmatrix}
i \sigma _k  & 0 \\
0 & - i \sigma _k 
\end{pmatrix} .
\end{equation} 
Hence positive chirality spinors belong to the   $\operatorname{Cl} (3) $ representation generated by
\begin{equation}
\label{gamma3rep} 
\gamma _k =i \sigma _k .
\end{equation}

In order to find  solutions of the eigenspinor equation (\ref{eigenspinoreqminus})  we look for harmonic spinor on $E ^4 $, which induce eigenspinors on $H ^3 $ as described in Section \ref{soleigenspinoreq}. The problem of finding solutions of the Dirac equation on $E ^4 $ in the background of an instanton has been solved in \cite{jackiw:1977}. In Appendix \ref{jrcomputation} we recall some of the results of \cite{jackiw:1977} and summarise the calculations needed to convert a specific $E ^4 $ harmonic spinor, associated to the JNR pole at the origin, into a spinor $\nu \in E _\xi  ^- $. The result is 
\begin{equation}
\label{numain} 
\nu = \frac{i}{2} (\nu _+ \sigma _+ + \nu _- \sigma _- )+ i\nu _3 \sigma _3 ,
\end{equation} 
where $\sigma _\pm =\sigma _1 \pm i \sigma _2 $ and $\nu _+ , \nu _- , \nu _3 $ are the Lie algebra components of $\nu$, given by
\begin{align} 
\label{numainplus} 
\nu _+  &
=\frac{4 \rho ^{ 5/2 } }{f  (x ^2 + y ^2 + \rho ^2 )^2  } 
\begin{pmatrix}
C _\rho   \\
C _x + i C _y 
\end{pmatrix} ,\\
\label{numainminus} 
\nu _- &
= \frac{4  (x - i y)  \rho ^{ 3/2 }}{f (x ^2 + y ^2 + \rho ^2 )^2  }
 \begin{pmatrix}
- C _x + i C _y  \\
C _\rho  
\end{pmatrix} ,\\
\label{numain3} 
\nu  _3 &=
 - \frac{2 \rho ^{ 3/2 } }{f (x ^2 + y ^2 + \rho ^2 )^2 } 
  \left[ (x- i y  ) 
  \begin{pmatrix}
C _\rho\\
C_ x + i C _y 
\end{pmatrix}  + \rho \begin{pmatrix}
C _x - i C _y \\
-C _\rho 
\end{pmatrix}  \right] ,
\end{align} 
with
\begin{equation}
\begin{split}
C _x &=\frac{\partial _x f}{f} + \frac{2x}{(x ^2 + y ^2 + \rho ^2 )^2 }, \ 
C _y =\frac{\partial _y f}{f} + \frac{2 y }{(x ^2 + y ^2 + \rho ^2 )^2 },\ 
C _\rho  = \frac{\partial _\rho f}{f} + \frac{2\rho }{(x ^2 + y ^2 + \rho ^2 )^2 } .
\end{split} 
\end{equation} 
This solution is the one described by equations (\ref{e4muplus})--(\ref{h3nu}) in Appendix \ref{jrcomputation} for the choice $ p _k =(x _k , y _k , z _k , w _k )= (0,0,0,0 ) $ of JNR pole, with the JNR function $f$ given by (\ref{jnrrefconfig}).  By making a different choice for the pole $p _k $ we could get other two more independent solutions.  In fact, as already mentioned, the full results in \cite{jackiw:1977} are sufficient to produce four explicit linearly independent solutions of the eigenspinor equation, that is a basis of $E _\xi  ^- $. However it is computationally much simpler to take a different  route as described in Section \ref{sn:inversionaction} below.

\subsection{An eigenspinor basis}
\label{sn:inversionaction}
The map $\invfour$ in (\ref{fourdiminvdef}) is an isometry of $H ^3 \times S ^1 $.  So if $B $ is $\invfour$-invariant and $\tilde{\nu}$ is a harmonic spinor on $H^3\times S^1$ with weight $\pm\frac12$ then $\invfour(\tilde{\nu})$ is also a harmonic spinor with weight $\mp\frac12$.  (The weight changes sign because $\invfour_\ast \partial_\theta=-\partial_\theta$).  Using the correspondence between harmonic spinors and eigenspinors on $H^3$, this gives a map from solutions of (\ref{eigenspinoreqminus}) to solutions of (\ref{eigenspinoreqplus}), which we now describe in detail.

The action of $\invfour $ on a 4-spinor on $H ^3 \times S ^1 $ in a local frame is written
\begin{equation} 
\invfour (\tilde \nu )  = \invlift \invfour ^\ast (\tilde \nu  ), 
\end{equation} 
for $\invlift $  a  $ \operatorname{Spin} (4) $ lift of the $SO (4) $  action of  $\invfour $ on an orthonormal coframe. Let us work with the upper half-space model of $H ^3 $ so that $\invfour (x,y, \rho , \theta ) = (I (x,y, \rho ), - \theta )$ with $I: H ^3 \rightarrow H ^3 $ given by (\ref{inv3uhs}). We take the $H ^3 \times S ^1 $ orthonormal coframe 
\begin{equation}
\epsilon ^1 =\frac{ \mathrm{d} x}{\rho }, \quad \epsilon ^2 =\frac{\mathrm{d} y}{\rho }, \quad \epsilon ^3 =\frac{\mathrm{d} \rho }{\rho }, \quad \epsilon ^4 =\mathrm{d} \theta .
\end{equation} 
Then with respect to the splitting of 4-spinors on $H ^3 \times S ^1 $ into positive and negative chirality, the lift of $\invfour $ to $ \operatorname{Spin} (4) $ is   $\invlift =\pm ( q \oplus q^{-1}  )$ with
\begin{equation}
q =\frac{1}{\sqrt{ x ^2 + y ^2 + \rho ^2 }} \begin{pmatrix}
\rho  & -x + i y \\
x + i y & \rho 
\end{pmatrix}.
\end{equation} 
Our results do not depend on the choice of the lift and from now on we take  $Q =q \oplus q ^{-1} $.
If $\tilde \nu=e^{\pm\frac{i}{2}\theta}\nu $  is a positive chirality spinor on $H ^3 \times S ^1 $ we then have
\begin{equation}
\label{ispintrasform} 
\invfour (\tilde \nu ) =
\invfour \left( \nu \mathrm{e}^{\pm\frac{i}{2}\theta}
\right) =
 q I ^\ast (\nu)\mathrm{e}^{\mp\frac{i}{2}\theta}.
\end{equation} 
Comparing (\ref{ispintrasform})  with (\ref{euclspinorweightp}), (\ref{eigenspinoreqderiv}) we  see that if  $\nu \in E _\xi  ^\pm $  then $qI^\ast(\nu)  \in E _\xi  ^\mp $. We will denote by $\Ispin: S \rightarrow S $ the map  induced by $\invfour $ on $ H ^3 $ spinors,
\begin{equation}
\Ispin  (\nu ) := q  I ^\ast (\nu).
\end{equation} 

In order to get another solution of (\ref{eigenspinoreqminus})  we need to combine the action of $\Ispin$ with that of the quaternionic structure on $S$ to obtain a quaternionic structure  on $ E _\xi  ^\pm $,
\begin{equation}
 E _\xi  ^\pm \rightarrow E _\xi  ^\pm , \quad \nu \mapsto \Ispin  (\thickbar{ \nu }) .
 \end{equation}
Hence, for the solution $\nu$ of (\ref{eigenspinoreqminus}) given by (\ref{numain})--(\ref{numain3}), the spinor $\Ispin(\thickbar\nu )$ also satisfies (\ref{eigenspinoreqminus}).

In calculations, one needs to keep in mind that $(A , \Phi )$ given in \eqref{jnrrefconfig}, \eqref{phimonopole}, \eqref{gae} is inversion symmetric only up to gauge transformations, see (\ref{invsymmonopole}). Therefore $ \Ispin (\thickbar \nu )$ satisfies (\ref{eigenspinoreqminus}) with respect to the gauge transformed monopole $(\operatorname{Ad} _{ \gi } A  + \gi \mathrm{d} \gi ^{-1} , \operatorname{Ad} _{ \gi }  \Phi  )$. In order to get a solution of (\ref{eigenspinoreqminus})  with respect to $(A , \Phi )$ we need to take $\operatorname{Ad} _{ \gi ^{-1}  }\Ispin (\thickbar \nu ) $ instead.  To ease on the  notation, from now on we will  write  $\Ispin(\thickbar \nu )$ when we are really referring to $\operatorname{Ad} _{ \gi ^{-1}  }\Ispin (\thickbar \nu ) $. For the monopole corresponding to the JNR function (\ref{jnrrefconfig}), see also Appendix \ref{jnrmonopoles}, the gauge transformation is 
\begin{equation}
\begin{split} 
\gi ^{-1} =
C ^{-1} \Bigg[
2   \rho  x \Big[ &x^4   x _0 ^2+2 x^2   \left(x _0 ^2 \left(\rho   ^2+y^2\right)-1\right)+x _0 ^ 2 \left(\left(\rho   ^2+y^2\right)^2+1\right)+2   \left(\rho ^2+y^2\right)\Big]I _2\\ 
+2x y \Big[ &x^4x _0 ^2+2 x^2\left(x _0 ^2 \left(\rho^2+y^2\right)-1\right)+x _0 ^   2 \left(\left(\rho   ^2+y^2\right)^2+1\right)+2   \left(\rho ^2+y^2\right)\Big]i \sigma _1 \\ 
-  \Big[ &x^6 x _0 ^2+x^4 \left(x _0 ^2   \left(\rho   ^2-x _0 ^2+y^2\right)-1\right   )-x^2 \Big(x _0 ^2   \left(\left(\rho ^2+y^2\right)   \left(\rho ^2+2   x _0 ^2+y^2\right)-1\right)  -6\left(\rho^2+y^2\right)\Big)\\ &
   -\left(\rho
   ^2+y^2\right) \left(\rho
   ^2+x _0 ^2+y^2\right)
   \left(x _0 ^2 \left(\rho
   ^2+y^2\right)+1\right)\Big]i  \sigma _2\Bigg],
   \end{split} 
\end{equation} 
with
\begin{equation}
\begin{split} 
C =(\rho ^2+x^2+y^2) \Big[ &
\left(\rho^2+(x-x _0 )^2+y^2\right)  \left(\rho   ^2+(x+x _0 )^2+y^2\right) \cdot \\ &
\left(x _0 ^2 \left(\rho   ^2+x^2+y^2\right)-2 x   x _0 +1\right) 
\left(x _0 ^2 \left(\rho^2+x^2+y^2\right)+2 xx _0 +1\right) 
\Big]^{ 1/2 }.
\end{split} 
\end{equation} 

We only need two more independent solutions  of (\ref{eigenspinoreqminus}) which we can obtain by Clifford acting with the monopole curvature on two linearly independent Killing spinors \cite{franchetti:2024}. For $\psi\in K ^- $ non-trivial, $\{ \psi , \Ispin (\thickbar{ \psi  }) \} $ is a basis of $K ^- $.
It can be checked that, with respect to half-space coordinates,
\begin{equation}
\label{psichoice} 
\psi =\begin{pmatrix}
\rho  ^{ 1/2 } \\
\rho ^{ -1/2 }(x + i y)\end{pmatrix}, \qquad 
\Ispin(\thickbar\psi ) = \begin{pmatrix}
0\\
- \rho ^{ -1/2 }
\end{pmatrix} 
\end{equation} 
are Killing spinors  satisfying
\begin{equation}
\label{psiinvpsipair} 
(\psi , \Ispin ( \thickbar\psi ) ) =-1.
\end{equation} 
So $F\cdot\psi$, $F\cdot\Ispin(\thickbar \psi)$ are independent solutions of \eqref{eigenspinoreqminus}. In conclusion, we take our basis of $E _\xi  ^- $  to be $\{ \nu , \Ispin (\thickbar \nu ),  F \cdot \psi  , F \cdot \Ispin (\thickbar \psi  ) \} $ with $\nu$ given in (\ref{numain})--(\ref{numain3}). Because of (\ref{isotangentspinors}),  $T ^{ \mathbb{C}  }_\xi  M ^f $ is then spanned by the eight tangent vectors
\begin{equation}
\label{alleigenspinors} 
\begin{split} 
\nu \otimes  \psi ,& \quad \nu \otimes  \Ispin (\thickbar\psi ), \quad \Ispin (\thickbar\nu ) \otimes \psi , \quad \Ispin (\thickbar\nu ) \otimes \Ispin (\thickbar\psi),\\
F \cdot \psi  \otimes  \psi ,& \quad F \cdot \psi  \otimes  \Ispin (\thickbar\psi ), \quad F \cdot \Ispin(\thickbar\psi)  \otimes \psi , \quad F \cdot  \Ispin(\thickbar\psi) \otimes \Ispin(\thickbar\psi).
\end{split} 
\end{equation}

In closing this Section we note the following identity, useful in many computations,
\begin{equation}
\thickbar{ \Ispin (\nu ) } =\Ispin (\thickbar { \nu } ).
 \end{equation} 
 In fact, 
 \begin{equation}
\begin{split}
\Ispin (\thickbar \nu  ) 
=q \begin{pmatrix}
I ^\ast (\bar\nu  _2 ) \\
-I ^\ast ( \bar\nu  _1 )
\end{pmatrix} ,
\end{split}
\end{equation} 
\begin{equation}
\begin{split}
\thickbar{\Ispin (\nu  )}
= \epsilon \bar{q}  \begin{pmatrix}
I ^\ast (\bar\nu  _1 ) \\
I ^\ast (\bar \nu  _2 )
\end{pmatrix} 
= q \epsilon \begin{pmatrix}
I ^\ast (\bar \nu  _1 ) \\
I ^\ast (\bar\nu  _2 )
\end{pmatrix} 
=q  \begin{pmatrix}
I ^\ast (\bar\nu  _2 ) \\
- I ^\ast (\bar\nu  _1 )
\end{pmatrix} ,
\end{split} \end{equation}
having used
\begin{equation}
q \epsilon =\epsilon \bar{q} ,
\end{equation}  
since, for $q \in SU (2) $, $q ^T  \epsilon q  =\epsilon   $ and $ \bar{q} ^T =q ^{-1} $. By the same argument, the action of the quaternionic structure commutes with that of any $SU (2) $ element.

We also note that the action of $\Ispin $  on  $(a, \phi ) \in T ^{ \mathbb{C}  }_\xi  M ^f $  induced by the isomorphism (\ref{isotangentspinors}) is
\begin{equation}
\label{aphiI} 
(a, \phi ) \mapsto (I ^\ast a , - I ^\ast \phi ) 
\quad \Leftrightarrow \quad 
(\nu , \psi ) \mapsto (\Ispin (\nu) , \Ispin (\psi) ) .
\end{equation}
(It is not simply pullback by $I$ on both $a $ and $\phi$, because $I$ is orientation-reversing so does not lift to $ \operatorname{Spin } (3) $).
For $u =(a, \phi ) \in T _\xi  ^{ \mathbb{C}  }M ^f $ we write
\begin{equation}
\label{iontangentspace} 
I (u) := (I ^\ast a , - I ^\ast \phi ).
\end{equation}

\section{Computing the $L^2$ metric}
\label{l2metric} 
\subsection{A symmetry-generated basis}
In order to compute the $L ^2 $ bilinear form  (\ref{L2product}) we start by choosing a convenient basis for $T ^ \mathbb{C}  _\xi M ^f $. The isometry group  $SO (3,1 )$ of $H ^3 $ acts on $M ^u $ by isometries with codimension one orbits.  
In order to compute the metric we can work at any point along the $SO (3,1) $ orbits. We choose the points $\xi =(A , \Phi )$ corresponding to the JNR function (\ref{jnrrefconfig}).

The directions tangent to the orbits correspond to infinitesimal rotations or boosts of the monopole. The remaining tangent direction in $M ^u $  corresponds to an infinitesimal motion separating the monopoles. The framed space $M ^f $ is an $\mathbb{R}  $-bundle over $M ^u $ with the framed direction generated by the infinitesimal unframed gauge transformation  $(\mathrm{d} ^A \Phi , 0 )$. It seems natural to choose a basis of tangent vectors generated by these monopole transformations.

Let $ (X _i , Y _i )$ be generators of $\mathfrak{so }(3,1 ) $, 
\begin{equation} 
[ X _i , X _j ] =- \epsilon _{ ijk }X _k , \quad [ X _i , Y _j ] =- \epsilon _{ ijk }Y _k , \quad [ Y _i , Y _j ] =\epsilon _{ ijk }X _k .
\end{equation} 
We use the same notation for the associated Killing vector fields on $H ^3 $, which in upper half-space coordinates are given by
\begin{equation}
\label{kvfh3upper} 
\begin{split}
X _1 &=  y  Y _3  +  \frac{1}{2} ( 1 - |x|  ^2  )\partial _y  ,\\
X _2 &=  -xY _3   - \frac{1}{2} (1 - |x|  ^2  )\partial _x  ,\\
X _3 &=  x  \partial _y  - y \partial _x,\\
Y _1 &=  - x Y _3   + \frac{1}{2} (1 + |x| ^2  ) \partial _x  ,\\
Y _2 &=  - y  Y _3   + \frac{1}{2} (1  + |x| ^2  ) \partial _y ,\\
Y _3 &=  x  \partial _x + y \partial _y + \rho  \partial _ \rho ,
\end{split}
\end{equation} 
where $|x |^2 =x ^2 + y  ^2 + \rho ^2 $.
Further overloading the notation, if $V$ is a Killing vector field we will also denote by $V$ the solution of (\ref{linbog}),  (\ref{gaugefixingminus}) corresponding to an infinitesimal rotation or boost. More precisely, there is some infinitesimal complex framed gauge transformation $(\mathrm{d} ^A \Lambda , [ \Phi , \Lambda ])$ such that $V =(a_V + \mathrm{d} ^A \Lambda , \phi _V+ [ \Phi , \Lambda ] )$, with
\begin{equation}
\label{infsym} 
a_V =L _V A = \iota  _V F+ \mathrm{d} _A \iota  _V A, \qquad \phi _V=L _V \Phi  =\iota  _V \mathrm{d} _A \Phi+ [ \Phi , i _V A ]  . 
\end{equation}
This gives six independent vectors $ \{ X _1 , X _2 , X _3 , Y _1 , Y _2 , Y _3 \} $ in $ T_\xi M ^f $.  A seventh is given by
\begin{equation}
\label{y0def} 
Y _0 =(\mathrm{d} ^A \Phi , 0 ),
\end{equation}
which is tangent to the fibres of the fibration $M^f\to M^u$.  Note that $Y_0$ solves (\ref{linbog}),  (\ref{gaugefixingminus}) because $(A,\Phi)$ solves the Bogomolny equation $d^A\Phi=\ast F$.  We complete the basis for $T_\xi M^f$ by choosing
\begin{equation}
\label{X0def} 
X _0 = \left( \frac{\mathrm{d} A}{\mathrm{d} x _0 }   + \mathrm{d} ^A \Lambda ,\frac{\mathrm{d} \Phi}{\mathrm{d} x _0 }    + [ \Phi , \Lambda ] \right) 
\end{equation}
for some suitable unframed gauge transformation $\Lambda$.  This is a lift to the framed moduli space $M^f$ of a tangent vector to $M^u$ obtained by varying the parameter $x_0$ in the JNR monopole \eqref{refconfigjnr}.  There is a freedom in the choice of lift which corresponds to adding multiples of $Y_0$ to $X_0$ (or equivalently, changing the choice of $\Lambda$ in \eqref{X0def}).  We will fix this freedom by imposing
\begin{equation}
\label{prodx0y0} 
g (X _0 , Y _0 )=0,
\end{equation}  
where $g$ is the $L ^2 $ bilinear form  (\ref{L2product}).  It can be checked that the vector $X _0 ^\prime = (\mathrm{d}  A / \mathrm{d} x _0, \mathrm{d} \Phi / \mathrm{d} x _0 )$ already satisfies this condition, that is
\begin{equation}
\label{X0framed} 
g (X _0 ^\prime , Y _0 )=0.
\end{equation} 
Hence imposing (\ref{prodx0y0}) is equivalent to requiring that $\Lambda$ in (\ref{X0def}) is an infinitesimal framed gauge transformation.

\label{theactionofd2}
Our choice of monopole $\xi$ given in \eqref{jnrrefconfig} is invariant under the action of $D_2$ given in \eqref{D2onS4}.  With respect to upper half-space coordinates on $H ^3 $ the $D _2 $ generators become
\begin{equation}
\label{D2generatosh3} 
\begin{split}
R _1 (x,y, \rho ) &= \frac{1 }{x ^2 + y ^2 + \rho ^2 } (x,-y, \rho ),\\
R _2 (x,y, \rho ) &= \frac{1 }{x ^2 + y ^2 + \rho ^2 } (-x, y, \rho ),\\
R _3 (x,y, \rho ) &=  (-x,-y, \rho ).
\end{split}
\end{equation}

Since $\xi$ is $D_2$-invariant, these act on the tangent space $T_\xi M^f$.  The action on our basis is:
\begin{equation}
R_i(X_k)=\begin{cases} X_k & k=i\text{ or }0\\-X_k & \text{otherwise},\end{cases}\qquad
R_i(Y_k)=\begin{cases} Y_k & k=i\text{ or }0\\-Y_k & \text{otherwise.}\end{cases}
\end{equation}
Because of these symmetry properties,  the $L ^2 $ bilinear form (\ref{L2product}) on $T _\xi ^{ \mathbb{C}  } M ^f $  simplifies considerably, taking the form
\begin{equation}
\label{l2metricgeneralform} 
g = h^2 \mathrm{d} x _0  ^2 + a _i ^2  \sigma _i^2  + b   _i ^2  \tau  _i ^2 + 2c _i ^2  \sigma _i \tau _i  + k ^2 \mathrm{d} t ^2  ,
\end{equation} 
with
\begin{alignat}{3}
\label{l2products} 
a _i ^2 &=g (X _i, X _i ),& \quad 
b _i ^2 &=g (Y _i , Y _i ),& \quad 
c _i ^2 &=g (X _i , Y _i ), \\ \nonumber 
 h^2 &=g (X _0 , X _0 ),& \quad
k ^2 &= g (Y _0 , Y _0 ).
\end{alignat} 
Here $t$ is a global coordinate along the framed direction, $\sigma _i $, $\tau _i $ are the left-invariant 1-forms on $SO (3,1) $ dual to the Killing vector fields $X _i $, $Y _i $, $i =1,2,3 $,  and $x _0 $ is the parameter in (\ref{jnrrefconfig}). Because of the $SO (3,1 )$ isometric action, $x _0 $ provides  a radial coordinate transverse to the $SO (3,1 )$ orbits and $a _i $, $b _i $, $c _i $,  $h$, $k$ are functions of $x _0 $ only.

\subsection{Change of basis}
The symmetry-generated basis $\{X_k,Y_k\}$ is the most natural basis in which to write the metric $g$.  However calculating the metric in this basis is difficult, because doing so entails solving the gauge-fixing condition \eqref{gaugefixingminus}, which amounts to finding gauge transformations $\Lambda$ for each basis vector which solve the appropriate PDE.

We will circumvent this difficulty by relating the basis $\{X_k,Y_k\}$ to the basis \eqref{alleigenspinors}.  The latter basis automatically solves the gauge-fixing condition \eqref{gaugefixingminus} since it was obtained from solutions of the Dirac eigenspinor equation \eqref{eigenspinoreqminus} and Killing spinors.  We will then calculate the metric in the  basis \eqref{alleigenspinors} and finally rewrite it in terms of the symmetry-generated basis $\{X_k,Y_k\}$.

We will use the action of $D_2$ on the tangent space to help us find the relations between the two bases.  
If $\nu\in E _\xi  ^- $, $\psi \in K ^- $, $ R \in D _2 \subset SO (3)  $ and $ \tilde R $ is a lift of $R$ to  $\operatorname{Spin} (3) =SU (2) $  then
\begin{equation} 
\label{d2d2staraction} 
R (\nu \otimes \psi ^\ast ) = \tilde R (\nu ) \tilde R (\psi^\ast ).
\end{equation} 
Since $\operatorname{Spin} (3) \rightarrow SO (3) $ is a double covering, both $\pm \tilde R $ are lifts of $R$, but since the lift of $R$ appears twice on the right hand side of  (\ref{d2d2staraction}), the  action on $\nu \otimes \psi ^\ast $ does not depend on the lift choice.

By lifting the action induced by the $D _2 $ generators (\ref{D2generatosh3})    on the orthonormal coframe $ \{ \rho ^{-1} \mathrm{d} x, \rho ^{-1} \mathrm{d} y , \rho ^{-1} \mathrm{d} \rho \} $ one obtains the   $ \operatorname{Spin} (3)  $  elements $\pm \tilde R _i $ given by
\begin{equation}
\begin{split}
\tilde R _1 &= - \frac{1}{\sqrt{ x ^2 + y ^2 + \rho ^2 }} (y I _2  + x i \sigma _3 - \rho i \sigma _1 ),\\
\tilde R _2 &=  \frac{1}{\sqrt{ x ^2 + y ^2 + \rho ^2 }} (x I _2  + \rho  i \sigma _2 - y i \sigma _3 ),\\
\tilde R _3 &= i \sigma _3 .
\end{split} 
\end{equation} 
The group generated by $ \{ \tilde R _i \} $ is the order eight binary dihedral group $ \tilde D _2 $, which is described by the relations $ \tilde R _i ^2 =-1 $, $ \tilde R _i \tilde R _j =- \tilde R _j \tilde R _i $ for $i \neq j $. For $\psi$ given by (\ref{psichoice}) one then has e.g.
\begin{equation}
\begin{split}
\tilde R _1 (\Ispin (\thickbar\psi ) )
=  -\left(  \frac{y I _2  + x i \sigma _3 - \rho i \sigma _1 }{\sqrt{ x ^2 + y ^2 + \rho ^2 }} \right) R _1   ^{-1\ast}  ( \Ispin (\thickbar\psi ) )
= - \frac{i }{\sqrt{ \rho }} \begin{pmatrix} 
\rho \\ x + i y
\end{pmatrix} = - i \psi .
\end{split} 
\end{equation} 
Applying $\tilde R _1 $ to both sides we get  $\tilde R _1 (\psi ) = -i \Ispin   (\thickbar\psi ) $. Similarly one calculates
\begin{alignat}{2}
 \tilde R _1 (\psi ) &= -i \Ispin   (\thickbar\psi ), &\quad \quad \tilde R _1 ( \Ispin  (\thickbar\psi ) ) &= - i \psi ,\\
\tilde R _2 (\psi ) &=  \Ispin   (\thickbar\psi ), &\quad \tilde R _2 ( \Ispin  (\thickbar\psi ) ) &=-\psi ,\\
\tilde R _3 (\psi ) &= i  \psi, &\quad \tilde R _3 ( \Ispin  (\thickbar\psi ) ) &=- i\Ispin  (\thickbar \psi) .
\end{alignat} 
This gives the action of $\tilde{D}_2$ on our basis for $K$.

Next we calculate the action on our basis for $E_\xi^-$. In general, for $O \in SO (3) $, $\tilde O $ one of its lifts to $SU (2) $, 
\begin{equation}
\tilde O (F \cdot \psi ) = O ^{-1\ast }(F ) \cdot \tilde O (\psi ).
\end{equation} 
 Since the monopole field strength $F$ is invariant under $D _2 $, up to gauge transformations
 \begin{equation}
 \tilde R _k (F \cdot \psi ) =F \cdot  \tilde R _k( \psi) , \qquad 
  \tilde R _k (F \cdot \Ispin(\thickbar \psi) ) =F \cdot  \tilde R _k( \Ispin(\thickbar\psi)) 
\end{equation}
so  $F \cdot \psi $, $F \cdot \Ispin (\thickbar\psi ) $ transform under $\tilde D _2 $, up to gauge transformations, in the same way as $\psi $, $\Ispin (\thickbar \psi ) $ do.

For $\nu$ given by (\ref{numainplus})--(\ref{numain3}) it is not hard to see that
\begin{equation}
\tilde R _3 (\nu _\pm )= i \nu _\pm , \quad 
\tilde R _3 ( \nu _3)= -i \nu _3 .
 \end{equation}
Conjugation by $i\sigma _3 $ maps $ \sigma _\pm \mapsto - \sigma _\pm $, hence up to gauge transformations,
\begin{equation} 
\tilde R _3 (\nu)=  -i \nu .
\end{equation} 
Since $ \Ispin \circ \tilde R _k = \tilde R _k \circ \Ispin $,   $\thickbar { \tilde R _k (\nu ) } = \tilde R _k (\thickbar \nu )$, it follows
\begin{equation}
\tilde R _3 (\Ispin   (\thickbar\nu ))
= i \Ispin   (\thickbar\nu ).
\end{equation} 
A lengthier calculation shows that, up to gauge transformations,
\begin{equation}
\tilde R _2 (\nu )= \Ispin (\thickbar\nu ).
\end{equation} 
The action of $\tilde R _1 $ is then
\begin{equation}
\label{akshjdkhjadhjkads} 
\tilde R _1 (\nu ) = \tilde R _3 (\tilde R _2 (\nu )) =i \Ispin  (\thickbar\nu ).
\end{equation} 
Collecting all the results, the action of $\tilde{D}_2$ on $E_\xi^-$ is
 \begin{alignat}{2}
 \label{nurotrasf1} 
 \tilde R _1 (\nu )&= i\Ispin  (\thickbar \nu ),&\quad \quad   \tilde R _1 (\Ispin  (\thickbar \nu ) ) &=  i \nu ,\\
  \label{nurotrasf2} 
\tilde R _2 (\nu )&= \Ispin  (\thickbar \nu ),&\quad \tilde R _2 (\Ispin  (\thickbar \nu ) ) &=-\nu ,\\
 \label{nurotrasf3} 
\tilde R _3 (\nu )&= -i \nu, &\tilde R _3 (\Ispin  (\thickbar \nu ) ) &= i \Ispin  (\thickbar \nu ),\\
 \tilde R _1 (\psi ) &= -i \Ispin   (\thickbar \psi ), &\quad \quad \tilde R _1 ( \Ispin  (\thickbar \psi ) ) &= - i \psi ,\\
 \label{r2psi} 
\tilde R _2 (\psi ) &=  \Ispin   (\thickbar \psi ), &\quad \tilde R _2 ( \Ispin  (\thickbar \psi ) ) &=-\psi ,\\
 \label{r3psi} 
\tilde R _3 (\psi ) &= i  \psi, &\quad \tilde R _3 ( \Ispin  (\thickbar \psi ) ) &=- i\Ispin  (\thickbar  \psi) .
\end{alignat}

Let us define
\begin{alignat}{2}
 \label{allinv} 
w _0  &= \nu \otimes  \psi ^\ast +  \Ispin  (\thickbar\nu ) \otimes\Ispin  (\thickbar\psi )^\ast ,& \qquad &v _0 =F \cdot \psi \otimes\Ispin  (\thickbar\psi )^\ast - F \cdot \Ispin  (\thickbar\psi )\otimes\psi ^\ast ,\\
\label{inv3} 
w _3  & =\nu \otimes\psi ^\ast - \Ispin  (\thickbar\nu ) \otimes\Ispin  (\thickbar\psi )^\ast ,& \qquad &v _3 =F \cdot \psi\otimes \Ispin  (\thickbar\psi )^\ast +  F \cdot \Ispin  (\thickbar\psi )\otimes\psi ^\ast ,\\
\label{inv1} 
w _1  & =\nu \otimes\Ispin  (\thickbar\psi) ^\ast +  \Ispin  (\thickbar\nu ) \otimes\psi^\ast ,& \qquad &v _1 =F \cdot \psi \otimes\psi ^\ast -  F \cdot \Ispin  (\thickbar\psi )\otimes\Ispin  (\thickbar\psi) ^\ast ,\\
\label{inv2} 
w _2  &= \nu \otimes\Ispin  (\thickbar\psi) ^\ast -  \Ispin  (\thickbar\nu ) \otimes\psi^\ast ,& \qquad &v _2 =F \cdot \psi \otimes\psi ^\ast +   F \cdot \Ispin  (\thickbar\psi )\otimes\Ispin  (\thickbar\psi)^\ast .
 \end{alignat}
This gives a basis for $T_\xi^\mathbb{C} M^f$ that is closely related to the basis \eqref{alleigenspinors}.  Using (\ref{nurotrasf1})--(\ref{r3psi})  one finds that, up to gauge transformations, $v _0 $, $w _0 $ are $D _2 $ invariant  while $ v _i $, $w _i $ are invariant under $R _i $ and change sign under $ R _j $ for $j \neq i $,  $j \in  \{ 1,2,3 \} $. Moreover one can check that
\begin{alignat}{4}
 \label{ibarw} 
 I (\bar w _1 ) &= - w _1,& \quad  I (\bar w _3 ) &= - w _3,& \quad 
  I (\bar w _2 ) &=   w _2,& \quad I ( \bar w _0 ) &= w _0,\\
  \label{ibarv} 
 I (\bar v _1 ) &= +  v _1,& \quad  I (\bar v _3 ) &= +  v _3,& \quad 
  I (\bar v _2 ) &= -  v _2,& \quad I ( \bar v _0 ) &= - v _0.
 \end{alignat} 
We recall that the moduli space   tangent vector associated to $v _i $, $w _i $  are determined by (\ref{spinorstotv}).

We are now ready to rewrite the basis $\{ v _i , w _i \} $ in terms of  $ \{ X _i , Y _i \} $. By comparing the way they transform  under  $D _2 $ we see that we must have
 \begin{equation}
w  _k = A _k  X _k  + B_k Y _k , \qquad v _k = C _k X _k +  D _k Y _k
\end{equation} 
for some constants $A _k , B _k , C _k , D _k $. To find the value of these constants we make use of the  identity
 \begin{equation}
 \label{gividentity} 
- \frac{1}{2} L _V \left(  \operatorname{Tr} (\Phi ^2 )  \right) 
= - \operatorname{Tr} (\Phi \phi  _V ).
 \end{equation} 
where $V$ is a Killing vector field and $ \phi _V  $ is given by (\ref{infsym}).  The rhs of \eqref{gividentity} does not change when $\phi_V$ is replaced by $\phi_V+[\Phi,\Lambda]$, so it is a well-defined function on $T_\xi^\mathbb{C}M^f$.

 For example, consider the tangent vector $w _3 $. Applying (\ref{spinorstotv}) to (\ref{inv3}) we calculate the corresponding value of $\phi  $, which we then plug in the rhs of (\ref{gividentity}). We calculate the lhs of (\ref{gividentity}) for $V = A _3 X _3 + B _3 Y _3 $  using the expressions for $X _3 $, $Y _3 $ given in (\ref{kvfh3upper}). Equating the two quantities then fixes the value of  $A _3 $, $B _3 $. The Higgs field $\Phi$ in (\ref{gividentity}) is the one of our chosen reference configuration, given by (\ref{phimonopole}) for $f$ the JNR function (\ref{jnrrefconfig}).  Performing these computations  one finds
\begin{align} 
\label{w1xy} 
w _1 &= -\left(  \frac{ 2( x _0^2 + 1 )^2 }{x _0 ^2 ( x _0 ^2 + 3 )} \right) Y _1 +2 i \left( \frac{  x _0 ^2 + 1 }{x _0 ^2} \right) X _1 ,\\
\label{w2xy} 
w _2 &=2 i \left( \frac{ x _0 ^2 -1}{x _0 ^2 } \right) Y _2 +  \left( \frac{ 2( x _0 ^2-1 )^2 }{x _0 ^2 (x _0 ^2 -3 )}\right) X _2,\\
\label{w3xy} 
w _3 &= \left( \frac{ x _0 ^4-1 }{ x _0 ^4 } \right) Y _3 ,\\
\label{v1xy} 
 v _1 &
= -  ( X _1  +  i Y _1  ) ,\\
\label{v2xy} 
  v _2 &
= -i  \left(X _2  +  i Y _2   \right) ,\\
\label{v3xy} 
 v _3 &
 = -  \left( X _3 +  i Y _3  \right)  .
 \end{align} 
It can be checked that the rhs of \eqref{w2xy} is finite as $x_0\to\sqrt{3}$, even though it involves a divergent quantity $1/(x_0^2-3)$. We can invert (\ref{w1xy})--(\ref{v3xy})   to get
\begin{align}
\label{X1vw} 
X _1 &= \frac{1}{2} (x _0 ^2 + 1 )v _1 -  \frac{i x _0 ^2 }{4} \left( \frac{ x _0 ^2 + 3 }{x _0 ^2 + 1} \right) w _1 ,\\
X _2 &= - \frac{ i}{2} (x _0 ^2 -3 ) v _2  + \frac{x _0 ^2 }{4} \left( \frac{x _0 ^2 -3 }{x _0 ^2 -1}  \right) w _2 ,\\
X _3 &=- v _3  -i\left( \frac{x _0 ^4 }{ x _0 ^4 -1 } \right) w _3 ,\\
Y _1 &=  \frac{ i}{2} (x _0 ^2 + 3 ) v _1 + \frac{x _0 ^2 }{4} \left( \frac{x _0 ^2 + 3}{x _0 ^2 + 1 } \right) w _1 ,\\
Y _2 &= \frac{1}{2}  (x _0 ^2 -1 ) v _2 + \frac{ix _0 ^2 }{4} \left( \frac{x _0 ^2 - 3}{x _0 ^2 - 1 } \right) w _2,\\
\label{Y3vw} 
Y _3 &=\left( \frac{ x _0 ^4 }{x _0 ^4 -1 } \right) w _3 .
\end{align}
Note that $X _2 $ vanishes for the equilateral configuration $x _0 =\sqrt{ 3 }$ of the poles, as it should since in that case rotation about the $x _2 $-axis is a symmetry of the monopole.

Equation (\ref{gividentity}) is sufficient  to establish the equalities (\ref{w1xy})--(\ref{v3xy}).  We double checked our results making use of a similar identity,
\begin{equation}
\label{givaidentity} 
- \frac{1}{2} L _V \operatorname{Tr} (F ^2 ) = - \operatorname{Tr}  \langle\mathrm{d} ^A a _{ V }, F  \rangle,
\end{equation} 
where $a _V $ is given by  (\ref{infsym}) and the monopole field strength $F$ by (\ref{monopfieldstrength}) for $f$ the JNR function (\ref{jnrrefconfig}).

For the tangent vector $w _0 $ we can apply (\ref{gividentity})  to $V =\mathrm{d} / \mathrm{d} x _0 $. Doing so however does not fully determine  $w _0 $  in terms of $X _0 $, $Y _0 $, giving
\begin{align}
\label{w0con} 
w _0 &=\left( \frac{x _0 ^4 -1}{x _0 ^3 } \right) X _0 + B _0 Y _0  .
\end{align} 
The constant $B _0 $ will be determined in Section \ref{thel2metric}. Note that  (\ref{givaidentity}) cannot be used to fix it since if $a=\mathrm{d} ^A \Phi $ then  
$\operatorname{Tr}  \langle \mathrm{d} ^A  a   , F \rangle =\operatorname{Tr}  \langle [ F , \Phi ]  , F \rangle  =0 $. 

For the tangent vector $v _0  $ calculating $(a, \phi )$ using equation (\ref{spinorstotv}) gives   $\phi =0 $ and, comparing  with (\ref{monstarF}),  $ a =- \mathrm{d} ^A \Phi $. By (\ref{y0def})  we thus have
\begin{equation}
\label{v0con} 
v _0 = - Y _0  =(- \mathrm{d} ^A \Phi ,0 ).
\end{equation} 

Behaviour under inversion offers a consistency check of our results. Our reference monopole $\xi =( A , \Phi )$ is inversion  symmetric so the map $ T ^ \mathbb{C}_\xi M ^f \rightarrow T^ \mathbb{C} _\xi M ^f $, $u \mapsto I (\bar u )$ is a real linear involution with eigenvalues $\pm1$, such that  $X_i$ have eigenvalue $ + 1 $ and $Y_i$ have eigenvalue $-1$ for $i =0,1,2,3 $.  Using (\ref{ibarw})--(\ref{ibarv}) one can verify that (\ref{X1vw})--(\ref{Y3vw}) transform correctly.  As for (\ref{w0con}) we find the constraint $B _0 \in i \,  \mathbb{R}$, which is indeed satisfied, see equation (\ref{B0value}) below.

Looking at the definitions of $\{ v _i \} $ in (\ref{allinv}) and at equations (\ref{v1xy})--(\ref{v3xy}), (\ref{v0con}) we see that suitable linear combinations of $F \cdot \psi $, $F \cdot  \Ispin(\thickbar \psi ) $ give the tangent vectors  $Y _0 $, $X_k - i Y _k $. This fact  had  already been derived, with different methods, in \cite{franchetti:2024}  and provides a further check of our results --  in \cite{franchetti:2024} the tangent vectors obtained are constant multiples of $Y _0 $ and $X _k + i Y _k $, the different sign   compared to (\ref{v1xy})--(\ref{v3xy})  being due to the opposite choice  $\gamma _k =- i \sigma _k $ made for the Clifford algebra representation.

\subsection{The $L^2$ metric}
\label{thel2metric}
We finally come to the computation of the $L ^2 $ form $g$.
We recall  \cite{franchetti:2024} that if $ u _i  =\nu_i \otimes \psi _i ^\ast \in T ^ \mathbb{C}_\xi  M ^f $ then the $L ^2$ form (\ref{L2product}) becomes
\begin{equation}
g (u _1 , u _2 )
= -\frac{1}{4}  (\psi _1 , \psi _2 ) \int _{ H ^3 } \operatorname{Tr} \left(  \nu _1 , \nu _2  \right)\operatorname{vol} ,
\end{equation} 
having used the fact that the symplectic pairing of Killing spinors is constant. In upper half-space coordinates the $H ^3 $ volume element is $\operatorname{vol} = \rho ^{ -3 } \mathrm{d} x \wedge \mathrm{d} y \wedge \mathrm{d} \rho $. Let us write
\begin{equation}
\omega (\nu _1 , \nu _2 ):= - \frac{1}{2} \int _{ H ^3 }\operatorname{Tr} (\nu _1 , \nu _2 ) \operatorname{vol} .
\end{equation} 
Note that from (\ref{L2product}), (\ref{barontangentspace}), (\ref{iontangentspace}) and the fact that $I$ is an isometry of $H ^3 $ follows that
\begin{equation}
\label{gbarinv} 
g (I (\bar u ), I (\bar v) )= \overline{g (u,v ) }.
\end{equation}

 From (\ref{allinv})--(\ref{inv2}) we obtain
\begin{equation}
\label{vvprod} 
  g (v _0 , v _0 )=g (v _2 , v _2 )=-g (v _1 , v _1 ) = -g (v _3 , v _3 ) .
\end{equation} 
Since the monopole energy is 
 \begin{equation} 
 \label{monenergy} 
E = - \frac{1}{2} \int _{ H ^3 }\operatorname{Tr} ( \mathrm{d} ^A \Phi \wedge * \mathrm{d} ^A \Phi )= 2 \pi p n,
 \end{equation} 
and in our case $p =\tfrac{1}{2} $, $n =2 $, it follows
\begin{equation}
\label{v0v0prod}
g (v _0 , v _0 )= g ((-\mathrm{d} ^A \Phi ,0 ), (-\mathrm{d} ^A \Phi , 0 ) )=E = 2 \pi .
\end{equation}

The symmetry properties of $\{ v _i , w _i \} $ and the invariance of $g$ under the $D _2 $ action imply 
\begin{equation}
\label{D2relations} 
g (w _i , w _j ) =   g (v _i , v _j ) =g (v _i , w _j )=0\quad  \text{for }  i \neq j, \  i,j \in \{ 0,1,2,3 \} .
\end{equation} 
Substituting (\ref{allinv})--(\ref{inv2}) one finds that the relations $g (v _i , v _j )=g (w _i , w _j )=0  $ are satisfied identically, while the conditions $g (v _i , w _j )=0 $ give 
\begin{align*} 
0 &=g (w _1 , v _2 )=g (w _2 , v _1 )=-g (w _0 , v _3 )=- g (w _3 , v _0 )
=\frac{1}{2}  \omega (\nu , F \cdot \psi )  - \frac{1}{2} \omega   (\Ispin  (\thickbar \nu) , F \cdot \Ispin  (\thickbar \psi)),\\
0 &=g (w _0 , v _1 ) =g (w _2 , v _3 )=- g (w _3 , v _2 )=- g (w _1 , v _0 ) 
=\frac{1}{2}  \omega (\nu , F \cdot \Ispin  (\thickbar \psi) )  + \frac{1}{2}  \omega  (\Ispin  (\thickbar \nu) , F \cdot \psi),\\
0 &=g (w _1 , v _3 )  =g (w _3 , v _1 )=- g (w _0 , v _2 ) = - g (w _2 , v _0 )
=\frac{1}{2}  \omega (\nu , F \cdot \Ispin  (\thickbar \psi) )  - \frac{1}{2} \omega (\Ispin  (\thickbar \nu) , F \cdot \psi).
\end{align*} 
It follows that
\begin{align}
\label{relD3R1} 
\omega (\nu , F \cdot \psi ) &=\omega  (\Ispin  (\thickbar\nu ), F \cdot \Ispin  (\thickbar\psi ) ),\\
\label{relD3R3} 
\omega (\nu , F \cdot \Ispin  (\thickbar\psi ) ) &=\omega  (\Ispin  (\thickbar\nu ), F \cdot \psi ) =0.
\end{align}   

The remaining products are
\begin{align}
\label{pr1} 
g (w _3 , w _3 )&=g (w _1 , w _1 )=-g (w _0 , w _0 )=-g (w _2 , w _2 )=\omega  (\nu , \Ispin  (\thickbar\nu ) ),\\
\label{pr2} 
g (w _1 , v _1 ) &=g (w _2 , v _2 ) = -g (w _0 , v _0 ) =-g (w _3 , v _3 ) = \omega (\nu , F \cdot \psi ),
\end{align} 
having also used (\ref{relD3R1}). By  (\ref{gbarinv}),  (\ref{ibarw}) 
\begin{align}
g (w _i , w _i )=\omega  (\nu , \Ispin  (\thickbar\nu ) )&\in \mathbb{R}  , \qquad 
g (w _i , v _i ) = \omega (\nu , F \cdot \psi )\in i \mathbb{R}  .
\end{align} 
In summary the non-zero products are
\begin{align}
\label{prodvv} 
g (v _0, v _0 ) &=g (v _2 , v _2 )=-g (v _1 , v _1 ) = -g (v _3 , v _3 )  = 2 \pi ,\\
\label{prodww} 
g (w _1 , w _1 )&=g (w _3 , w _3 )=-g (w _0 , w _0 )=-g (w _2 , w _2 ) = \omega (\nu , \Ispin (\thickbar\nu ) ) \in \mathbb{R}  ,\\
\label{prodwv} 
g (w _1 , v _1 ) &=g (w _2 , v _2 ) = -g (w _0 , v _0 ) =-g (w _3 , v _3 )= \omega  (\nu , F \cdot \psi ) \in i \mathbb{R}  ,
\end{align} 
leaving only two spinor pairings to be computed,
\begin{align}
\label{intnuibnu} 
\omega (\nu , \Ispin (\thickbar\nu ) ) &
= - \frac{1}{2}  \int _{ H ^3 } \operatorname{Tr} (\nu , \Ispin (\thickbar\nu ) ) \operatorname{vol} ,\\
\label{intnunukill} 
\omega (\nu , F \cdot \psi  ) &
= - \frac{1}{2}  \int _{ H ^3 } \operatorname{Tr} (\nu , F \cdot \psi  ) \operatorname{vol} .
\end{align}

We have been able to evaluate both integrals \eqref{intnuibnu}, \eqref{intnunukill} explicitly.  The computation is not straightforward, and we present our method in detail in Appendix \ref{integrals}. The result is
\begin{alignat}{3}
\label{intnuibnufinal}
\omega (\nu ,\mathscr I (\thickbar \nu )) & 
= &  c _1 \left[ E(m) - E\left(\varphi  ,\,m\right)\right] & &
+ c _2\left[ K(m) - F\left(\varphi  , \,m\right) \right] &,\\
\label{intnunukillfinal}
\omega (\nu , F \cdot \psi )& 
= & k _1 
\left[ E(m) - E\left(\varphi  ,\,m\right)\right] & &
+ k _2
\left[ K(m) - F\left(\varphi  , \,m\right) \right] & + 2\pi i\, \left( \frac{x _0 ^4 -1 }{x _0 ^4 }\right),
\end{alignat}
where $K (m )$, $F( \varphi , m ) $ are the complete and incomplete elliptic integrals of the first kind, $E(m) $, $E(\varphi , m )$ are the complete and incomplete elliptic integrals of the second kind, and
\begin{alignat}{2}
\varphi &=\arcsin \left( \sqrt{ 1-\frac{4}{(x_0^2+1)^2}} \right), & \quad  
m &=\frac{(x_0^2-3)(x_0^2 + 1)^3}{(x_0^2 + 3)(x_0^2-1)^3},\\
k _1 &=- 2 \pi i  \frac{(x _0 ^2 -1 )(x _0 ^4   + 3  )}{x _0 ^4 (x _0 ^2 -3) }  \sqrt{ \frac{ x _0 ^2 -1 }{ x _0 ^2 + 3 }} , &\quad 
k _2 &=  \frac{16\pi i}{x _0 ^2 (x _0 ^2 -3 )}  \ \,  \sqrt{ \frac{x _0 ^2 -1 }{ x _0 ^2 + 3 } },\\
c _1 &=-\frac{8\pi\,(x_0^2-1)^2(x_0^2+1)}{x_0^4(x_0^2-3)}
\sqrt{\frac{x_0^2-1}{x_0^2+3}}, &\quad 
c _2 &=\frac{32\pi\,(x_0^2+1)}{x_0^4(x_0^2-3)}
\sqrt{\frac{x_0^2-1}{x_0^2+3}}.
\end{alignat}
Interestingly, as shown in Appendix \ref{integrals}, the elliptic curve that appears in these elliptic integrals is biholomorphic to the spectral curve of the monopole.  It would be interesting to find an explanation for this surprising observation.

\begin{figure}[htbp]
\begin{center}
\includegraphics[scale=1]{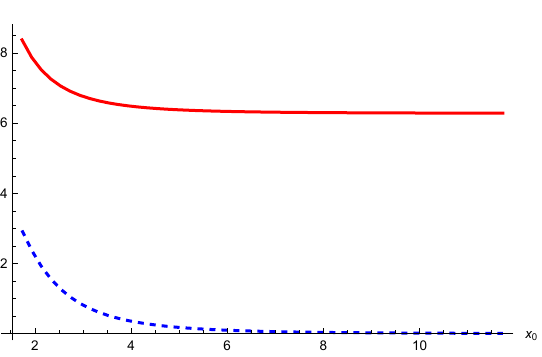}
\caption{The red solid, blue dashed lines show the value of $-i\omega (\nu , F \cdot \psi ) $, $\omega ( \nu , \Ispin (\thickbar \nu ) )$  as functions of $x _0 $. }
\label{integralplots}
\end{center}
\end{figure}

Figure \ref{integralplots} shows the plots of (\ref{intnuibnufinal}), (\ref{intnunukillfinal}). The relations (\ref{prodww}), (\ref{prodwv}) are satisfied as expected. The function  $\omega (\nu , \Ispin (\thickbar\nu ) )$ is non-zero real positive, decreasing with $x _0 $. The function $\omega ( \nu , F \cdot \psi )$ is  non-zero purely imaginary, with $| \omega (\nu , F \cdot \psi ) | = -i \omega (\nu , F \cdot \psi )$ decreasing with $x _0 $. The asymptotic expressions for large $x _0 $ are
\begin{align}
\label{intnuibnuasy} 
\omega (\nu ,\mathscr I (\thickbar \nu ))  &
\simeq \frac{16 \pi }{ x _0 ^4 } \left( 2 \log x _0 -1 \right) ,\\
\label{intnunukillasy}
\omega (\nu , F \cdot \psi ) &
\simeq 2 \pi  i + \frac{2 \pi  i}{x _0 ^4 } \left(  8 \log  x _0 -3  \right) .
\end{align} 

We can now fix the value of the constant $ B _0 $ in (\ref{w0con}). By (\ref{X0framed}), \eqref{v0con} \eqref{prodvv} and \eqref{prodwv}, taking the $L ^2 $ product of (\ref{w0con}) with $ Y _0 $ gives
\begin{equation} 
\label{B0value} 
 B _0 =\frac{ g ( w _0 , Y _0 )}{g(Y_0,Y_0) }=\frac{ -g ( w _0 , v _0 )}{g(v_0,v_0) }= \frac{ 1}{ 2 \pi }  \omega (\nu , F \cdot \psi ) \in i \, \mathbb{R}  .
 \end{equation}
  Therefore
\begin{equation}
\label{w0x0temp} 
w _0 =\left( \frac{x _0 ^4 -1}{x _0 ^3 } \right) X _0 +  \frac{1}{2 \pi }\omega  (\nu , F \cdot \psi ) Y_0 ,
\end{equation} 
or equivalently
\begin{equation}
X _0 =\left( \frac{x _0 ^3}{x _0 ^4 -1 } \right) \left( w _0 +  \frac{1}{2 \pi }\omega (\nu , F \cdot \psi ) v_0\right).
\end{equation}

We can now write down the products (\ref{l2products}) which determine the $L ^2 $ form (\ref{l2metricgeneralform}) on $T_\xi^\mathbb{C}M ^f  $. As expected,  the only non-zero products are  $g (X _i , X _i )$, $ g (Y _i , Y _i )$, $ g (X _i , Y _i )$. From (\ref{X1vw})--(\ref{Y3vw}) and  (\ref{prodvv})--(\ref{prodwv}) we finally get
\begin{align}
\label{gx0y0} 
g (X _0 , X _0 )&
= \left( \frac{x _0 ^3 }{x _0 ^4 -1 } \right) ^2 \left( \frac{1}{2 \pi } |\omega (\nu , F \cdot \psi )|^2 -\omega  (\nu , \Ispin (\thickbar\nu ) )\right)  ,\\
g (X _1 , X _1 ) &
=- \frac{ \pi }{2}  (x _0 ^2 + 1) ^2
- \frac{x _0 ^4 }{16} \left( \frac{x _0 ^2 + 3}{ x _0 ^2 + 1} \right)^2\omega  (\nu , \Ispin  (\thickbar\nu ) )
-  \frac{ i}{4}x _0 ^2 (x _0 ^2 + 3)\omega  (\nu , F \cdot \psi ),\\
g (X _2 , X _2 ) &
=- \frac{ \pi}{2} (x _0 ^2 -3) ^2 -  \frac{x _0 ^4 }{16} \left( \frac{x _0 ^2 - 3}{ x _0 ^2 - 1 } \right) ^2\omega  (\nu , \Ispin  (\thickbar\nu ) )
-\frac{ i}{4}  x _0 ^2 \left( \frac{(x _0 ^2 - 3 ) ^2  }{  x _0 ^2 -1 } \right)  \omega (\nu , F \cdot \psi ),\\
g (X _3 , X _3 ) &
= - 2 \pi -\left(  \frac{ x _0 ^4 }{x _0 ^4 -1  } \right)^2  \omega (\nu  ,\Ispin (\thickbar\nu)  )- \left( \frac{2 ix _0 ^4 }{x _0 ^4 -1 }\right) \omega  (\nu , F \cdot \psi ),\\
g (Y_0 , Y _0 )&=2 \pi ,\\
g (Y _1 , Y _1 ) &
= \frac{  \pi}{2}  \left(x _0 ^2 + 3 \right) ^2 +  \frac{x _0 ^4 }{16} \left( \frac{x _0 ^2 + 3}{ x _0 ^2 + 1} \right) ^2 \omega  (\nu , \Ispin  (\thickbar\nu ) )
+ \frac{ i}{4} x _0 ^2 \left( \frac{(x _0 ^2 + 3)^2 }{ x _0 ^2 + 1 } \right) \omega (\nu , F \cdot \psi ),\\
g (Y _2 , Y _2 ) &
=\frac{  \pi}{2}  \left( x _0 ^2 -1 \right) ^2 +  \frac{x _0 ^4 }{16} \left( \frac{x _0 ^2 - 3}{ x _0 ^2 - 1 } \right) ^2 \omega  (\nu , \Ispin  (\thickbar\nu ) )
+ \frac{ i}{4}  x _0 ^2 \left( x _0 ^2 - 3 \right) \omega  (\nu , F \cdot \psi ),\\
g (Y _3 , Y _3 ) &
=\left(  \frac{ x _0 ^4 }{x _0 ^4 -1  } \right)^2 \omega  (\nu  ,\Ispin (\thickbar\nu)  ).
\end{align} 
The expressions for the mixed products $g (X _i , Y _i )$, $i =1,2,3 $, follow from (\ref{v1xy})--(\ref{v3xy})   and (\ref{prodvv}),
\begin{align}
g (X _0 , Y _0 ) &= 0,\\
g(X _1, Y _1 )&
= -\frac{ i}{2} \pi    (x _0 ^2 + 1 ) (x _0 ^2 + 3 ) - \frac{ix _0 ^4 }{16} \left( \frac{x _0 ^2 + 3}{x _0 ^2 + 1} \right) ^2\omega  (\nu , \Ispin (\thickbar\nu ))
+\frac{1}{4}   \left( \frac{ x _0 ^2 (x _0 ^2 + 2 ) (x _0 ^2 + 3)}{x _0 ^2 + 1 } \right) \omega (\nu , F \cdot \psi ),\\
g(X _2, Y _2 )&
= -\frac{ i}{2} \pi  (x _0 ^2 -1 ) (x _0 ^2 -3 )
- \frac{ix _0 ^4 }{16} \left( \frac{ x _0 ^2 - 3 }{x _0 ^2 - 1 } \right)^2  \omega (\nu , \Ispin (\thickbar\nu ))
+\frac{1}{4}   \left( \frac{x _0 ^2 (x _0 ^2 -2 )(x _0 ^2 - 3 ) }{x _0 ^2 -1 } \right)\omega  (\nu , F \cdot \psi )  ,\\
\label{gx3y3} 
g (X _3 , Y _3 )  &= -i \left( \frac{x _0 ^4 }{x _0 ^4 -1 } \right) ^2 \omega (\nu , \Ispin (\thickbar\nu )) +  \left(\frac{   x _0 ^4 }{x _0 ^4 -1} \right) \omega (\nu , F \cdot \psi ).
\end{align} 
As we can see from (\ref{l2metricgeneralform}), (\ref{l2products}), the $L ^2 $ form on the moduli space of charge 2, mass $\tfrac{1}{2}$ hyperbolic moduli is fully determined by (\ref{gx0y0})--(\ref{gx3y3}). As expected \cite{franchetti:2024}, the products $g (X _i , X _i ) $, $g (Y _i , Y _i )$ of vectors on the tangent space to inversion symmetric monopoles and on its normal bundle are real. The mixed products $g (X _i , Y _i )$ are purely imaginary.

\section{Asymptotics}
\label{asymptotics}
\subsection{Comparison with $H^3\times H^3$}
Having computed the $L^2$ metric, we now analyse its asymptotics.  As $x_0\to\infty$, a 2-monopole resembles a well-separated pair of 1-monopoles.  Each of these 1-monopoles is described by its location in $H^3$ and a phase parameter.  The natural geometry to expect as $x_0\to\infty$ is that of a torus bundle (or $\mathbb{R}\times S^1$-bundle) over $(H^3\times H^3)/\mathbb{Z}_2$.  Here $\mathbb{Z}_2$ swaps the two 1-monopoles, reflecting the fact that 1-monopoles are indistinguishable.  To investigate whether this is true of the $L^2$ metric, we first write the metric on $H^3\times H^3$ in a suitable form.

The isometry group $SO(3,1)$ of $H^3$ acts diagonally on $(H^3\times H^3)/\mathbb{Z}  _2 $ with cohomogeneity 1.  The orbits are represented by points
\begin{equation}
((0,0,\rho_0), (0,0,1/\rho_0))\in H^3\times H^3,\quad 0<\rho_0\leq1,
\end{equation}
written in upper half-space coordinates $(x,y,\rho)$.  The tangent space at any such point is spanned by $\mathrm{d}/\mathrm{d}\rho_0$ and by the tangent vectors that generate the action of $SO(3,1)$.  The latter take the form
\begin{equation}
\big(X_i(0,0,\rho_0),X_i(0,0,1/\rho_0),\,
Y_i(0,0,\rho_0),Y_i(0,0,1/\rho_0)\big)\in T_{(0,0,\rho_0)}H^3\times T_{(0,0,1/\rho_0)}H^3\,,
\end{equation}
where $X_i,Y_i$ are the vector fields written explicitly in \eqref{kvfh3upper}.  We will abuse notation by referring to these tangent vectors as $X_i,Y_i$.  In this basis, the metric on $H^3\times H^3$ takes the form
\begin{equation}\label{H3H3 metric}
\frac{2}{\rho_0^2}\mathrm{d}\rho_0^2 + \frac{1}{2}\left(\rho_0-\frac{1}{\rho_0}\right)^2(\sigma_1^2+\sigma_2^2)
+ \frac{1}{2}\left(\rho_0+\frac{1}{\rho_0}\right)^2(\tau_1^2+\tau_2^2) + 2\tau_3^2.
\end{equation}
Here, as in \eqref{l2metricgeneralform}, $\sigma_i,\tau_i$ are left-invariant 1-forms dual to $X_i,Y_i$.  The small $\rho_0$ asymptotics of (\ref{H3H3 metric}) is
\begin{equation}\label{H3H3 asymptotic metric}
\frac{2}{\rho_0^2}\mathrm{d}\rho_0^2 + \frac{1}{2\rho_0^2}(\sigma_1^2+\sigma_2^2+\tau_1^2+\tau_2^2) + 2\tau_3^2.
\end{equation}

In order to compare this with the metric on $M$, we need to identify the locations of the well-separated 1-monopoles.  We do so using spectral curves.  The spectral curve of a monopole written in JNR form \eqref{rrc} is \cite{bolognesi:2015}
\begin{equation}
0=\sum_{j=1}^3\lambda_j^2\prod_{k\neq j}(z-x_k-iy_k)(1+(x_k-iy_k)w).
\end{equation}
For the reference monopole with parameters given in \eqref{refconfigjnr}, this is
\begin{equation}
\label{2-monopole curve}
0=(1-x_0^4)wz+x_0^2(w^2+z^2)-(1+w^2z^2).
\end{equation}
This form of the curve agrees up to a $90^\circ$ rotation with the curve found in \cite{sutcliffe:2022}.  More precisely, on making the substitutions
\begin{equation}
x_0^2=\frac{3+a}{1-a},\quad w\mapsto iw,\quad z\mapsto iz
\end{equation}
it agrees with the curve (3.5) in \cite{sutcliffe:2022}.

We wish to compare this with the spectral curves of well-separated 1-monopoles.  The spectral curve of a 1-monopole located at $(0,0,\rho_0)$ is $z=\rho_0^2w$.  So the union of the spectral curves of a pair of 1-monopoles at $(0,0,\rho_0^{\pm1})$ is
\begin{equation}
0=(z-\rho_0^2w)(w-\rho_0^2z)\approx wz-\rho_0^2(w^2+z^2)+O(\rho_0^4)
\end{equation}
This agrees with \eqref{2-monopole curve} up to $O(\rho_0^4)$ on making the identification
\begin{equation}
\rho_0=\frac{1}{x_0}.
\end{equation}
So the 1-monopole locations are asymptotically $(0,0,x_0^{\pm1})$ as $x_0\to\infty$.

Having identified the asymptotic 1-monopole locations, we can compare the asymptotic $L^2$ metric with the metric on $H^3\times H^3$.  The asymptotic $L^2$ metric obtained by taking the leading coefficients in \eqref{gx0y0}--\eqref{gx3y3}
\begin{equation}\label{asymptotic metric}
\h =\frac{2\pi}{\rho_0^2}\mathrm{d}\rho_0^2+\frac{\pi}{2\rho_0^2}(\sigma_1^2+\sigma_2^2+\tau_1^2+\tau_2^2)+2\pi \tau_3^2+2\pi(\sigma_3+i\tau_3)^2 + 2\pi\mathrm{d}t^2.
\end{equation}
The first six terms are equal to $\pi$ times the metric \eqref{H3H3 asymptotic metric} on $H^3\times H^3$.  The factor of $\pi$ is equal to the energy \eqref{monenergy} of a monopole with charge $n=1$ and mass parameter $p=\frac12$.  The coefficient of the final term $\mathrm{d}t^2$ is equal to the total energy of the 2-monopole, and this is true not just asymptotically but on the whole moduli space (see \eqref{v0v0prod} and \cite{franchetti:2024}).  So these parts of the metric are natural and agree with expectations.

The seventh term $2\pi(\sigma_3+i\tau_3)^2 $ in the asymptotic metric \eqref{asymptotic metric} is more surprising.  It can be interpreted in terms of the relative phase of the 1-monopoles.  The vector field $X_3$ generates an action of $S^1$ which fixes the 1-monopole positions and changes their relative phase.  Acting with $\exp(\pi X_3)$ is the same as acting with $R_3\in D_2$ and is a symmetry of the monopole.  In the asymptotic metric \eqref{asymptotic metric} this circle has circumference $\sqrt{2\pi^3}$.

Globally, these $S^1$ orbits form a circle bundle over the subset of $H^3\times H^3$ corresponding to well-separated monopoles.  The metric \eqref{asymptotic metric} is a Kaluza-Klein metric on this circle bundle, because we can write $\sigma_3+i\tau_3=\mathrm{d}\theta + A$ in terms of a locally-defined coordinate $\theta$ on the $S^1$ fibres and a locally-defined gauge potential $A$ on $H^3\times H^3$.  The field strength (or curvature) of the gauge potential is
\begin{equation}
\mathrm{d}A = \mathrm{d}(\sigma_3+i\tau_3)=-(\sigma_1+i\tau_1)\wedge(\sigma_2+i\tau_2).
\end{equation}
The fact that this is complex (not real) is a reflection of the fact that the $L^2$ metric is complex.  It would be interesting to try to understand why this field strength is complex by adapting Gibbons--Manton's analysis of the asymptotic metrics of Euclidean monopoles \cite{gibbons:1995}.

\subsection{Asymptotics of the centred moduli space}
It is also interesting to look at the asymptotic metric on the moduli space $M ^I $  (\ref{modulispacecentred}) of centred 2-monopoles.  This is the submanifold of $M$ obtained by taking the orbit of our reference monopole \eqref{jnrrefconfig} under the action of the subgroup $SO(3)\subset SO(3,1)$.  It consists of all 2-monopoles invariant under the inversion $I$.  The asymptotic metric on this submanifold was predicted in \cite{franchetti:2023aa} to be a circle bundle over $H^3$.

The asymptotic $L^2$ metric on the centred moduli space is obtained by restricting the asymptotic $L^2$ metric \eqref{asymptotic metric}:
\begin{equation}\label{centred asymptotic metric}
2\pi \left( \frac{1}{\rho_0^2}\mathrm{d}\rho_0^2+\frac{1}{4\rho_0^2}(\sigma_1^2+\sigma_2^2)+\sigma_3^2 \right) .
\end{equation}
The metric on $H^3$ is
\begin{equation}
\frac{1}{\rho_0^2}\mathrm{d}\rho_0^2 + \frac{1}{4}\left(\rho_0-\frac{1}{\rho_0}\right)^2(\sigma_1^2+\sigma_2^2)=\frac{1}{\rho_0^2}\mathrm{d}\rho_0^2 + \left(\frac{1}{4\rho_0^2}+O(1)\right)^2(\sigma_1^2+\sigma_2^2).
\end{equation}
This can be obtained similarly obtained by restricting the metric \eqref{H3H3 metric} on $H^3\times H^3$ to the inversion-symmetric submanifold, and dividing by 2.  Dividing by 2 ensures that the metric has scalar curvature $-1$.  Alternatively, one can calculate this metric directly by choosing representatives $(0,0,\rho_0)\in H^3$ of the orbits of the $SO(3)$ action.  To leading order, the $L^2$ metric \eqref{centred asymptotic metric} is that of a circle bundle over $H^3$, in agreement with the predictions of \cite{franchetti:2023aa}.

It is interesting to compare the asymptotic metric \eqref{asymptotic metric} with Hitchin's self-dual Einstein metrics on moduli spaces of centred hyperbolic 2-monopoles \cite{hitchin:1996}.  The moduli space of centred hyperbolic 2-monopoles is diffeomorphic to $S^4\setminus\mathbb{RP}^2$.  Hitchin constructed a family of self-dual Einstein metrics on this moduli space using the twistor picture of hyperbolic monopoles. These metrics depend on an integer parameter $k \geq 3 $ and, for $k>3 $, have conical defects with angle $2 \pi / (k-2 )$ along $\mathbb{RP}^2\subset S^4$.

Multiplying the asymptotic metric \eqref{centred asymptotic metric} by the conformal factor $\rho_0^2$ gives a metric with a conical defect along the $\mathbb{RP}^2 \subset S ^4 $.  So the $L^2$ metric will be asymptotically conformal to Hitchin's metric if and only if the cone angle is the same as Hitchin's metric.

We can compute the cone angle by computing the radius of curvature $r_{H^3}$ of the hyperbolic metric and the radius $r_{S^1}$ of the circle fibres in \eqref{centred asymptotic metric}.  The cone angle will be equal to $2\pi r_{S^1}/r_{H^3}$ (and the metric will be continuous if $r_{S^1}=r_{H^3}$).  As we have already seen, the circle fibres are generated by the action of $X_3$, and $\exp(\theta X_3)$ acts nontrivially for $0<\theta<\pi$ and trivially for $\theta=\pi$.  So the restriction of the metric \eqref{centred asymptotic metric} to the circle fibres is $2\pi$ times that of a circle of circumference $\pi$, in other words, it is the metric on a circle of radius $r_{S^1}=\sqrt{{\pi}/{2}}$.  The asymptotic metric on $H^3$ in \eqref{asymptotic metric} is $2\pi$ times the metric with sectional curvature $-1$, so its radius of curvature is $r_{H^3}=\sqrt{2\pi}$.  
So the conical defect angle
\begin{equation}
2\pi\frac{r_{S^1}}{r_{H^3}}=\pi.
\end{equation}
Hitchin's metric on the moduli space of centred hyperbolic monopoles 2-monopoles with mass parameter $\frac12$ is the $k =6 $ member of the family and has conical defect angle $\pi/2$.  So the $L^2$ metric is not asymptotically conformal to Hitchin's metric for the same mass parameter.  Amongst all of the metrics constructed in \cite{hitchin:1996}, the only one with conical defect angle $\pi$ is the Fubini--Study metric on $\mathbb{CP}^2/\mathbb{Z}_2$, obtained for $k =4 $, where $\mathbb{Z}_2$ acts by complex conjugation and fixes $\mathbb{RP}^2\subset\mathbb{CP}^2$.  So the $L^2$ metric is asymptotically conformal to the Fubini--Study metric, but we have checked that it is not conformal to the Fubini-Study metric globally.

\subsection{The Gibbons--Warnick metric}
Gibbons--Warnick have computed a metric for well-separated  charge $n$ monopoles in which all but one of the monopoles have electric charge 0 and fixed locations.  In this section we compare our asymptotic metric with theirs in the case $n=2$.

In order to do this, we need to fix the location of one monopole and set its electric charge to zero.  In the moduli space \eqref{framed moduli space} monopoles are indistinguishable.  In order to control the electric charge and position of one monopole we need to treat them as distinguishable, so we work instead on
\begin{equation}\label{distinguishable moduli space}
M=(C,\infty)\times SO(3,1)/\Z_2\times \R.
\end{equation}
Here $C$ is a large real number and $\Z_2\subset D_2$ is generated by the rotation $R_3$ that fixes the monopole locations.  The other two rotations in $D_2$, $R_1$ and $R_2$, swap the monopole locations.   The space \eqref{distinguishable moduli space} is a double cover of an open subset of the moduli space \eqref{framed moduli space} in which the monopoles are distinguishable.

Now let us set the electric charge of one monopole to zero.  We can do this by quotienting by the action generated by  the vector field $Z_2$ that changes the phase of the second monopole.  We claim that the two vector fields that change the phases of each of the monopoles are
\begin{equation}
\label{phase vector fields} 
Z_1=\frac{1}{2}\left(\frac{d}{dt}+X_3\right),\quad Z_2=\frac{1}{2}\left(\frac{d}{dt}-X_3\right).
\end{equation}
Here ${X}_3$ is the left-invariant vector field on $SO(3,1)/\Z_2$ which is dual to the left-invariant 1-form $\sigma_3$.

To prove (\ref{phase vector fields})  we first note that 
\begin{equation}
Z_1+Z_2=\frac{d}{dt}
\end{equation}
since $\frac{d}{dt}$ corresponds to the tangent vector $(a,\phi)=(d^A\Phi,0)$, which changes the phases of all monopoles.  On the other hand, $Z_1-Z_2$ must be a constant multiple of $X_3$ since $X_3$ is the only vector field on the moduli space that changes neither the total phase nor the positions of the monopoles.  We argue that in fact
\begin{equation}
\label{Zdiff} 
Z_1-Z_2=X_3.
\end{equation}

To prove (\ref{Zdiff}) we  consider the double-cover of the unframed moduli space of well-separated  distinguishable monopoles. It is   the quotient of \eqref{distinguishable moduli space} by the action generated by $\frac{d}{dt}$, i.e.~
\begin{equation}
\label{asydisti} 
(C,\infty)\times SO(3,1)/\Z_2.
\end{equation}
This space  inherits the actions generated by $Z_1$ and $Z_2$.

Consider the action of $\exp(2\pi Z_1)$ on (\ref{asydisti}).  Near the second monopole it acts trivially.  Near the first monopole it acts as a gauge transformation $g=\exp(2\pi\Phi)$.  Far from both monopoles, this gauge transformation is approximately constant and equal $-1$, which acts trivially on the monopole fields.  We extend the  gauge transformation to the whole of $H^3$ in such a way that it equals $-1$ near the second monopole.  Then the action of $\exp(2\pi Z_1)$ is equivalent to a gauge transformation with $g$, so acts trivially on the unframed moduli space.  In contrast, for $0<\theta<2\pi$ the action of $\exp(\theta Z_1)$ is not equivalent to a gauge transformation. Since $Z_1+Z_2$ acts trivially on the unframed moduli space, it follows that $\exp(\theta(Z_1-Z_2))$ acts trivially if and only if $\theta\in\pi\Z$.  On the other hand, $\exp(\theta X_3)$ acts trivially on \eqref{distinguishable moduli space} if and only if $\theta\in\pi\mathbb{Z}$ because the generator of $\Z_2$ is $R_3=\exp(\pi X_3)$.  Thus $Z_1-Z_2=\pm X_3$, and up to relabelling, $Z_1-Z_2=X_3$.

We now return to calculating the metric that describes motion where the second monopole has electric charge 0.  Since the asymptotic metric \eqref{asymptotic metric} is a Kaluza--Klein metric, setting the electric charge to 0 is equivalent to taking the quotient metric for the action generated by $Z_2$\footnote{Here is a justification written in local coordinates.  A Kaluza--Klein metric can be written in the form $g_{\mu\nu}dx^\mu dx^\nu + (d\theta+A_\mu(x) dx^\mu)^2$.  The geodesic equations are $\dot{q}=0$ and $\frac{d}{dt}(g_{\mu\nu}\dot x^\nu) -\frac{1}{2} \partial _\mu g _{ \nu \rho \dot{x} ^\nu \dot{x} ^\rho }+q(\partial_\nu A_\mu-\partial_\mu A_\nu)\dot x^\nu=0$, in which $q=\dot\theta+A_\mu\dot x^\mu$ is the Noether charge associated with $\frac{\partial}{\partial\theta}$.  Setting $q=0$ reduces these to $\frac{d}{dt}(g_{\mu\nu}\dot x^\nu)-\frac{1}{2} \partial _\mu g _{ \nu \rho }\dot{x} ^\nu \dot{x} ^\rho=0 $, which are the geodesic equations of the quotient metric $g_{\mu\nu}dx^\mu dx^\nu$.}.  The quotient of the moduli space \eqref{distinguishable moduli space} is
\begin{equation}\label{reduced distinguishable moduli space}
\bar{M}=(C,\infty)\times SO(3,1)/\Z_2.
\end{equation}
The quotient map is $\pi:M\to\bar{M}$, $\pi:(x_0,g\Z_2,t)\mapsto (x_0,g\exp(tX_3)\Z_2)$.  Then $\pi^\ast \sigma_3=\sigma_3+dt$, $\pi^\ast\sigma_i=\sigma_i$ for $i=1,2$, and $\pi^\ast\tau_i=\tau_i$ for $i=1,2,3$.
The quotient metric is
\begin{equation}\label{reduced distinguishable metric}
\hb=\frac{2\pi}{\rho_0^2}\mathrm{d}\rho_0^2+\frac{\pi}{2\rho_0^2}(\sigma_1^2+\sigma_2^2+\tau_1^2+\tau_2^2)+2\pi \tau_3^2+\pi(\sigma_3+i\tau_3)^2.
\end{equation}
This clearly satisfies $h=\pi^\ast\hb+\pi(\sigma_3+i\tau_3-dt)^2$, which confirms that $\hb$ is the quotient metric.

Now let us fix the position of one monopole.  The moduli space $\bar{M}$ in \eqref{reduced distinguishable moduli space} fibres over an open subset of $H^3\times H^3$, and the coordinates on the two copies of $H^3\times H^3$ represent the positions of the two monopoles.  We will fix the position of the second monopole to be $(x_2,y_2,\rho_2)=(0,0,1)$.  This constraint determines a four-dimensional submanifold $N$ of the moduli space \eqref{reduced distinguishable moduli space}.  We wish to calculate the restriction of the metric \eqref{reduced distinguishable metric} to this subspace.

The subgroup $SO(3)\subset SO(3,1)$ that fixes $(0,0,1)\in H^3$ acts on $N$ with cohomogeneity one.  This means that $N$ is the union of $SO(3)$-orbits of a 1-parameter family of points $p\in\bar{M}$ corresponding to monopoles with positions
\begin{equation}\label{GW base point}
(x_1,y_1,\rho_1,x_2,y_2,\rho_2)=(0,0,\rho,0,0,1),\quad 0<\rho\ll 1.
\end{equation}
We will calculate the metric on $N$ using this action, writing it in terms of the left-invariant 1-forms of $SO(3)$ and the coordinate $\rho$.  To do this, we start with a basis for $T_pN$.  The natural basis to choose is
\begin{equation}
\frac{\partial}{\partial\rho},\,X_1(p),\,X_2(p),X_3(p),
\end{equation}
in which $X_1,X_2,X_3$ are the vector fields on $\bar{M}$ that generate the left action of $SO(3)\subset SO(3,1)$.
Note that, although the points \eqref{GW base point} have an axial symmetry generated by $X_3$, the tangent vector $X_3(p)\in T_p\bar{M}$ is non-zero because the associated monopole is not axially-symmetric.

We now calculate the metric on $T_pN$ for a fixed point $p\in N$ of the form \eqref{GW base point} with $\rho$ fixed to some constant value $a$.
The metric \eqref{reduced distinguishable metric} is written in terms of points in $\bar{M}$ with position coordinates $(0,0,\rho_0,0,0,1/\rho_0)\in H^3\times H^3$.  So to calculate the metric on $T_pN$ we need to first find an isometry $f\in SO(3,1)$ of the metric \eqref{reduced distinguishable metric} that maps our chosen point $p$ to a point of this form.  A suitable choice is $f=\exp(-\frac12\ln a Y_3)$.  This acts on $H^3\times H^3$ as follows: 
\begin{equation}
\label{isopfp} 
f:(x_1,y_1,\rho_1,x_2,y_2,\rho_2)\mapsto \frac{1}{\sqrt{ a}} (x_1,y_1,\rho_1,x_2,y_2,\rho_2).
\end{equation}
In particular,   it maps $(0,0, a,0,0,1)$ to $(0,0,\sqrt{ a},0,0,1/\sqrt{ a})$.
The metric on $T_pN$ is then the restriction to $T_pN$ of the metric \eqref{reduced distinguishable metric} pulled back to $T_p \bar M $. It can be calculated by taking the inner products via $\hb$ of the tangent vectors spanning $T _p N $ pushed forward to $T _{ f (p)} \bar M $.

The pushforward of $X_3(p)$ is given by $f_\ast X_3(p)=X_3(f(p))$, because $[X_3,Y_3]=0$.   To calculate the pushforward of the remaining basis vectors we employ local coordinates on $H^3\times H^3$.  Our basis vectors correspond to the following tangent vectors to $H^3\times H^3$ at $(0,0, a,0,0,1)$:
\begin{equation}
\begin{aligned}
\frac{\partial}{\partial \rho} &= \frac{\partial}{\partial\rho_1} \\
X_1(p) &= \frac{1}{2}(1- a^2) \frac{\partial}{\partial y_1} \\
X_2(p) &= \frac{1}{2}( a^2-1) \frac{\partial}{\partial x_1}.
\end{aligned}
\end{equation}
We have written these in terms of coordinate vector fields using the explicit formulae \eqref{kvfh3upper} for the Killing vector fields of $H^3$.
Similarly, a basis for the tangent space of $H^3\times H^3$ at $(0,0,\sqrt{ a},0,0,1/\sqrt{ a})$ is given by
\begin{equation}
\label{tvecfp} 
\begin{aligned}
X_1(f(p))&=\frac12\left(1- a\right)\frac{\partial}{\partial y_1}+\frac12\left(1-\frac{1}{ a}\right)\frac{\partial}{\partial y_2}\\
X_2(f(p))&=\frac12\left( a-1\right)\frac{\partial}{\partial x_1}+\frac12\left(\frac{1}{ a}-1\right)\frac{\partial}{\partial x_2}\\
Y_1(f(p))&=\frac12\left( a+1\right)\frac{\partial}{\partial x_1}+\frac12\left(\frac{1}{ a}+1\right)\frac{\partial}{\partial x_2}\\
Y_2(f(p))&=\frac12\left( a+1\right)\frac{\partial}{\partial y_1}+\frac12\left(\frac{1}{ a}+1\right)\frac{\partial}{\partial y_2}\\
Y_3(f(p))&=\sqrt{ a}\frac{\partial}{\partial\rho_1}+\frac{1}{\sqrt{ a}}\frac{\partial}{\partial\rho_2}\\
\frac{\partial}{\partial \rho_0}(f(p))&
= \left. \frac{\partial \rho _1  }{\partial \rho _0 } \frac{\partial }{\partial \rho _1 } +\frac{\partial \rho _2 }{\partial \rho _0 } \frac{\partial }{\partial \rho _2} \right | _{ f (p) }
=\frac{\partial}{\partial \rho_1}-\frac{1}{a}\frac{\partial}{\partial \rho_2},
\end{aligned}
\end{equation}
since $f (p) = (0,0 , \rho _0 , 0,0, 1/ \rho _0 )$ with $\rho _0 =\sqrt{  a }$.
A short calculation shows that
\begin{equation}\label{pushforward}
\begin{aligned}
f_\ast \frac{\partial}{\partial \rho}&=\frac{1}{\sqrt{ a}}\frac{\partial}{\partial\rho_1}=\frac{1}{2\sqrt{ a}}\frac{\partial}{\partial \rho_0}(f(p))+\frac{1}{2 a}Y_3(f(p)),\\
f_\ast X_1(p)&=\frac{1- a^2}{2\sqrt{ a}}\frac{\partial}{\partial y_1}=\frac{1- a}{2\sqrt{ a}}Y_2(f(p))+\frac{1+ a}{2\sqrt{ a}}X_1(f(p)),\\
f_\ast X_2(p)&=\frac{ a^2-1}{2\sqrt{ a}}\frac{\partial}{\partial x_1}=\frac{ a-1}{2\sqrt{ a}}Y_1(f(p))+\frac{ a+1}{2\sqrt{ a}}X_2(f(p)),
\end{aligned}
\end{equation}
having used (\ref{tvecfp}) to express $\partial _{ x _1 }$, $\partial _{ y _1 }$, $\partial _{ \rho _1 }$ in terms of the vector fields $X _i $, $Y _i $.
These formulae allow us to calculate the components of the metric \eqref{reduced distinguishable metric} at $p$.  For example,
\begin{equation}
\begin{split} 
f ^\ast \hb \left(X _1 (p) , X _1 (p)   \right) &
=\hb\left(\frac{1- a}{2\sqrt{ a}}Y_2(f(p))+\frac{1+ a}{2\sqrt{ a}}X_1(f(p)),\frac{1- a}{2\sqrt{ a}}Y_2(f(p))+\frac{1+ a}{2\sqrt{ a}}X_1(f(p)) \right) \\ &
= \frac{(1- a)^2 }{4 a}  \hb(Y_2(f(p)),Y_2(f(p)))+ \frac{(1+ a)^2 }{4 a}  \hb(X_1(f(p)),X_1(f(p)))\\ &
= \frac{\pi }{2 \rho _0 ^2 }\frac{ (1 +  a ^2 )}{2  a } 
= \frac{\pi  (1 +  a ^2 )}{4  a ^2  } .
\end{split} 
\end{equation}
To calculate the metric at a general point of the form \eqref{GW base point}, rather than a fixed point $p$ with $\rho=a$, we just substitute $a =\rho $ into these formulae.  The result of the calculation is that,  to leading order, the asymptotic metric on $N$, obtained for $\rho \ll 1 $, is given by
\begin{equation}\label{GW metric}
{\pi}\frac{d\rho^2}{\rho^2}+\frac{\pi}{4\rho^2}\left(\bar{\sigma}_1^2+\bar{\sigma}_2^2\right)+\pi\left(\bar{\sigma}_3+\frac{id\rho}{2\rho}\right)^2.
\end{equation}
Here $\bar{\sigma}_1,\bar{\sigma}_2,\bar{\sigma}_3$ are left-invariant 1-forms on $SO(3)$ that are dual to $X_1,X_2,X_3$.  The bar in our notation distinguishes them from the left-invariant 1-forms $\sigma_i$ on $SO(3,1)$.  The metric \eqref{GW metric} describes the motion of well-separated 2-monopoles in which one of the monopoles has fixed position and electric charge 0.

We recognise the first three terms in the metric \eqref{GW metric} as the asymptotic metric on $H^3$.  The complex metric \eqref{GW metric} takes the form of a Kaluza--Klein metric on a circle bundle over a neighbourhood of the boundary in $H^3$.  In local coordinates, the metric takes the form $\pi g_{\mu\nu}dx^\mu dx^\nu+\pi(d\theta+A_\mu dx^\mu)^2$, in which $g_{\mu\nu}dx^\mu dx^\nu$ is the real metric on $H^3$ and $A=A_\mu dx^\mu$ is a locally-defined complex 1-form.  The geodesic equations are $\dot{q}=0$ and $\frac{d}{dt}(g_{\mu\nu}\dot x^\nu) -\frac{1}{2} \partial _\mu g _{ \nu \rho }\dot{x} ^\nu \dot{x} ^\rho +qF_{\nu\mu}\dot x^\nu=0$, in which $q=\dot\theta+A_\mu\dot x^\mu$ is the Noether charge associated with $\frac{\partial}{\partial\theta}$ and
\begin{equation}\label{GW field strength}
F = dA = d\left(\bar{\sigma}_3+\frac{id\rho}{2\rho}\right)=-\bar{\sigma}_1\wedge\bar{\sigma_2}.
\end{equation}
These are the equations of motion for a particle with electric charge $q$ in $H^3$ moving in the background of the real field strength tensor $F$.  Physically, $F$ is the field strength tensor of a Dirac magnetic monopole.  If the conserved electric charge $q$ is chosen to be real then the dynamics in $H^3$ will be real.  Hence the metric \eqref{GW metric}, while complex, induces real dynamics.

Gibbons--Warnick have also calculated a metric that describes the same situation of a monopole moving in the background of another monopole with fixed position and electric charge 0.  They derived this metric by treating the monopoles as point sources in classical electromagnetism.  Their metric was a hyperbolic Taub--NUT metric, with asymptotic form
\begin{equation}
{B}\frac{d\rho^2}{\rho^2}+\frac{B}{4\rho^2}\left(\bar{\sigma}_1^2+\bar{\sigma}_2^2\right)+C\bar{\sigma}_3^2
\end{equation}
for real constants $B,C$.  This is again a Kaluza--Klein metric and is clearly different from our asymptotic metric.  However, the associated dynamics on $H^3$ is the same: it is motion of a charged particle in the presence of the magnetic field \eqref{GW field strength}.  So while our metric is different from Gibbons--Warnick's, it induces the same dynamics for the position of the moving monopole.

\section{Conclusions}
\label{conclusions}
We have calculated the $L^2$ metric on the moduli space of hyperbolic monopoles in the particular case of charge 2 and mass parameter $\frac12$.  Our calculation shows that the metric is complex, not real, addressing a question raised in \cite{franchetti:2024}.

Our calculation reveals a simple asymptotic form for the metric: it is a Kaluza--Klein metric over $H^3\times H^3$.  Asymptotic metrics for euclidean monopoles are well-understood and can be derived from an interpretation of monopoles as point particles \cite{gibbons:1995}.  It would be interesting to develop a similar understanding of the asymptotic metric for hyperbolic monopoles.  So far this has been attempted only in two specific cases: where all but one of the monopoles are frozen \cite{gibbons:2007}, and for centred charge 2 monopoles \cite{franchetti:2023aa}.  The dynamics predicted by our metric is consistent with both of these studies.

The methods introduced in this paper could be used to calculate $L^2$ metrics for other charges and mass parameters.  In particular, they could be used to calculate the metric for charge $3$ and mass parameter $\frac12$, because in this case the moduli space is parametrised by the JNR ansatz.  Very little is known about moduli space metrics for 3-monopoles, either in euclidean space or in hyperbolic space.  So calculating the metric in this case would be interesting.

Our 2-monopole metric is written explicitly in terms of elliptic integrals whose underlying elliptic curve is biholomorphic to the monopole spectral curve.  This is unlikely to be a coincidence, and suggests that there is a way to express the $L^2$ metric for hyperbolic monopoles of any charge in terms of twistor theory.  It would be both interesting and useful to find such an expression.

Apart from the $L^2$ metric, another (possibly complex) metric for hyperbolic monopoles has been introduced by Nash using twistor theory \cite{nash:2007}.  Nash's metric has some structural similarities with the $L^2$ metric.  Specifically, both metrics are plurihermitian, meaning that they are compatible with an associated pluricomplex structure.  These similarities suggest that the $L^2$ and Nash metrics might be equivalent.  So it would be interesting to calculate the Nash metric and compare it with the $L^2$ metric calculated in this paper.  This will be done in another paper which is in preparation.

The gauge-fixing condition used to define the $L^2$ metric is motivated by supersymmetry.  This gauge-fixing condition arises as part of the BPS equations in a supersymmetric gauge theory on hyperbolic space, whose BPS states are monopoles.  So there should be some kind of supersymmetric quantum mechanics on the moduli space of hyperbolic monopoles that involves the $L^2$ metric.  It would be interesting to formulate this quantum mechanical theory precisely and use it to study the quantum dynamics of hyperbolic monopoles.

\appendix
\section{JNR monopoles}
\label{jnrmonopoles} 
To obtain the  monopole data corresponding to the JNR instanton  (\ref{cijnransatz})  we need to switch to polar coordinates $ \rho \mathrm{e} ^{ i \theta } =z + i w $ and to a gauge where $L _{ \partial _\theta }B=0$. Then $B =A - \Phi \mathrm{d} \theta $ with $A (\partial _\theta )=0 $. The required gauge change is 
\begin{equation}
B  \mapsto  g B g ^{-1} + g  \mathrm{d} g^{-1}, \quad g =\exp \left(  \frac{i}{2}  \theta \sigma _3 \right) .
\end{equation} 
One obtains
\begin{align}
\label{phimonopole} 
 \Phi &
=\frac{i \z}{2\rrc} \left( \partial _x \rrc \, \sigma _1  + \partial _y \rrc \, \sigma _2  +  \partial _\z \rrc  \sigma _3  \right) 
+  \frac{i }{2}  \sigma _3 ,
\\  
\label{gae}
 A  &
 =
-\frac{i}{2\rrc} \Big[
\left(  \partial _y \rrc \, \mathrm{d} \z -  \partial _\z \rrc  \mathrm{d} y \right)\sigma _1 +
\left(   \partial _\z \rrc  \mathrm{d} x - \partial _x \rrc \mathrm{d} \z \right) \sigma _2 + 
\left(  \partial _x \rrc \, \mathrm{d} y -  \partial _y \rrc  \mathrm{d} x \right) \sigma _3 
\Big ] .
\end{align} 
Using the fact that $f$ is harmonic, $ (\partial ^2 _x + \partial ^2 _y + \rho ^{-1} \partial _\rho + \partial ^2_ \rho )f =0 $, one can explicitly check that $(A , \Phi )$ satisfies $\mathrm{d} _A \Phi = * F $.

The field strength of (\ref{gae})  has components
 \begin{equation}
 \label{monopfieldstrength} 
 \begin{aligned}
F _{ xy }&
 = -\frac{i}{2} \Big[ \Big( (\partial ^2 _x + \partial ^2 _y ) \log \rr + (\partial _\rho  \log \rr )^2 \Big) \sigma _3  \\ \nonumber &
 + \Big(\partial _x \log \rr\partial _\rho  \log \rr  - \partial _{ x \rho  } \log \rr \Big) \sigma _1 
  + \Big(\partial _y \log \rr\partial _\rho  \log \rr  - \partial _{ y \rho  } \log \rr \Big) \sigma _2   \Big]\\
F _{ y\rho  }&
 = - \frac{i}{2} \Big[  \Big( ( \partial ^2 _\rho  + \partial ^2  _y ) \log \rr + (\partial _x \log \rr )^2  \Big ) \sigma _1 
 + \Big( \partial _x  \log \rr \partial _y \log \rr - \partial _{ xy } \log \rr \Big) \sigma _2  \\ \nonumber &
 + \Big( \partial _x \log \rr \partial _\rho  \log \rr - \partial _{x\rho  } \log \rr \Big) \sigma _3  \Big] \\
F _{ \rho x }&
 = -\frac{i}{2} \Big[  \Big( \partial _x \log \rr \partial _y \log \rr -\partial _{xy} \log \rr  \Big ) \sigma _1  
 + \Big( (\partial ^2  _\rho  + \partial ^2 _x )  \log \rr   +  (\partial _{ y } \log \rr ) ^2 \Big) \sigma _2  \\ \nonumber &
+ \Big( \partial _y  \log \rr \partial _\rho  \log \rr - \partial _{ y\rho  } \log \rr \Big)\sigma _3  \Big] ,
 \end{aligned}
 \end{equation} 
 while
 \begin{equation}
 \label{monstarF} 
 * F = \mathrm{d} _A \Phi = \rho (F _{ xy } \mathrm{d} \rho + F _{ y \rho } \mathrm{d} x + F _{ \rho x }\mathrm{d} y ).
 \end{equation}

\section{Harmonic spinors on $E^4$}
\label{jrcomputation} 
Let $B$ be an ASD   instanton with instanton number $n$ given by the JNR ansatz. Jackiw-Rebbi \cite{jackiw:1977} constructed $4n + 4 $ solutions to the Dirac equation $D ^B _E \nue =0$. The solutions  are not all linearly independent but rather satisfy four linear constraints, bringing the number of independent solutions to $4 n $, which is  the full dimension of the space of solutions.  The harmonic spinors constructed  in \cite{jackiw:1977} fall into two families, but we will only describe the simpler one since that will be sufficient for our purposes.  Let
\begin{equation}
f = \sum _{ k =1 } ^{ n + 1 } \frac{\lambda _k ^2 }{|p - p _k |^2 } 
\end{equation} 
be the JNR harmonic function on $E ^4 $ with $p =(x,y,z,w )$, $ p _k =(x _k , y _k  , z _k , w _k )$. Since $B$ is ASD, any spinor $\nue $  satisfying $D^B_E\nue=0$ has positive chirality, so we can  identify it with a 2-spinor.

A family of $2 (n + 1) $ independent solutions is then given by, writing $\nue= \nue   _a \, i\sigma _a $, $a =1, 2, 3 $, 
\begin{equation} 
\label{jrmksol} 
 \nue ^k   _a  = f \alpha   ^j \partial _j \left(  -i \gamma  _a \frac{  M  ^k  u}{ f ^2 } \right) , 
\end{equation} 
with, for $k =1 , \ldots , n + 1 $,
\begin{equation}
 M^k  = \frac{\lambda _k ^2 }{|p- p _k |^4 } \bar\alpha  ^j (p- p _k )_j 
 = \frac{ i \lambda _k ^2 }{|p- p _k |^4 }  \begin{pmatrix}
z- z _k -i (w- w _k)  &  x- x _k  -i (y - y _k ) \\
x- x _k  + i (y - y _k ) & - z + z _k -i (w - w _k )
\end{pmatrix} .
\end{equation} 
Here  $u$ is a constant 2-spinor, $  \alpha ^j = (-  \gamma  _1 , - \gamma _2 , - \gamma _3 , I _2 ) $, $ \bar{\alpha} ^j =(\gamma  _1, \gamma _2 , \gamma _3  , I _2 )$.  The $n + 1 $ choices of $ M _i $ and the two independent components of $u$ give $2 (n + 1 )$ solutions.  Write
 \begin{equation} 
 V = \begin{pmatrix}
V _1  \\V _2 
\end{pmatrix}= \begin{pmatrix}
z- z _k -i (w- w _k)  &  x- x _k  -i (y - y _k ) \\
x- x _k  + i (y - y _k ) & - z + z _k -i (w - w _k )
\end{pmatrix}  u ,
 \end{equation} 
\begin{equation}
\nue  
=\frac{i}{2} (\nue _- \sigma _+ + \nue _ +  \sigma _- ) + i \nue _3 \sigma _3 ,
\end{equation} 
where $\nue _\pm =\nue _1 \pm i \nue _2 $, $ \sigma _\pm =\sigma _1 \pm i \sigma _2 $.
Expanding (\ref{jrmksol})  one finds
\begin{align} 
\nue _+&
= - \frac{4 \lambda _k ^2 V _2 }{f |p-p _k |^4 }\begin{pmatrix}
 C _z + i C _w  \\
C _x + i C _y 
\end{pmatrix} ,\\
\nue_-&
= - \frac{4 \lambda _k ^2 V _1}{f |p-p _k |^4 }\begin{pmatrix}
C _x - i C _y  \\
- C _z + i C _w 
\end{pmatrix} ,\\
\nue _3 &=
 - \frac{2 \lambda _k ^2 }{f |p-p _k |^4 } \left[V _1 \begin{pmatrix}
C _z + i C _w  \\
C_ x + i C _y 
\end{pmatrix}  - V _2 \begin{pmatrix}
C _x - i C _y  \\
-C _z +  i C _w 
\end{pmatrix}  \right]  ,
\end{align} 
with
\begin{equation}
\label{cccdef} 
C _j 
= \frac{\partial _j f}{f} + \frac{2  (p - p _k)_j  }{|p - p _k | ^2 }.
\end{equation} 

We now assume that the JNR data is circle-invariant, $z _k =w _k =0 $, and introduce polar coordinates $\rho \mathrm{e}  ^{ i \theta }=z + i w$. 
Requiring the  spinor to have weight $ \tfrac{1}{2} $ under the circle action forces the first component of  $u$ to vanish, cutting down the number of solutions to $n + 1 $. Taking $=u =(0,1 )^T $ we get
\begin{equation} 
V =\begin{pmatrix}
x- x _k - i (y - y _k ) \\
- z + z _k - i (w-w _k )
\end{pmatrix} .
\end{equation} 
\begin{equation}
C_ z \pm  i C _w =\mathrm{e} ^{ \pm i \theta }C _\rho , \quad C _\rho =\frac{\partial _\rho f}{f} + \frac{2 \rho }{|p - p _k |^2 }.
\end{equation} 
The $SO (4) $ matrix relating the orthonormal coframes $\{ \mathrm{d} x , \mathrm{d} y , \mathrm{d} z ,  \mathrm{d} w \} $, $\{ \mathrm{d} x , \mathrm{d} y , \mathrm{d} \rho , \rho \mathrm{d} \theta \} $ lifts to the $\operatorname{Spin} (4) =SU (2) \times SU (2) $ matrix
\begin{equation}
\pm \operatorname{diag}\left( \mathrm{e} ^{  -\frac{ i}{2} \sigma _3 \theta } , \mathrm{e} ^{  \frac{ i}{2} \sigma _3 \theta } \right) .
\end{equation} 
Upon this change of frame $\nue $ thus becomes
\begin{align} 
\nue  _+ &
=\frac{4 \lambda _k ^2 \rho }{f |p-p _k |^4 } \mathrm{e} ^{3 i \theta /2}
\begin{pmatrix}
C _\rho   \\
C _x + i C _y 
\end{pmatrix} ,\\
\nue_- &
= -\frac{4 \lambda _k ^2 [x- x _k - i (y - y _k )]}{f |p-p _k |^4 }\mathrm{e} ^{ -i \theta /2} 
 \begin{pmatrix}
 C _x - i C _y  \\
-C _\rho  
\end{pmatrix} ,\\
\nue _3 &=
 - \frac{2 \lambda _k ^2}{f |p-p _k |^4 } \mathrm{e} ^{ i \theta /2} 
  \left[[x- x _k - i (y - y _k ) ] 
  \begin{pmatrix}
C _\rho\\
C_ x + i C _y 
\end{pmatrix}  + \rho \begin{pmatrix}
C _x - i C _y \\
-C _\rho 
\end{pmatrix}  \right] .
\end{align} 

We now switch to a circle-invariant gauge  conjugating by $\exp (- (i/2) \theta \sigma _3 )$, so that
$\sigma _3 \mapsto \sigma _3$, $\sigma _\pm \mapsto \mathrm{e} ^{ \mp  i \theta }\sigma _\pm$. We also drop the constant factor $\lambda _k ^2 $.
The gauge transformed spinor has then the form
$ \mathrm{e} ^{ i \theta /2 } \hat \nue$ with $\hat \nue  $  circle-invariant,
\begin{align} 
\label{e4muplus} 
\hat\nue  _+ &
=\frac{4  \rho }{f |p-p _k |^4 } 
\begin{pmatrix}
C _\rho   \\
C _x + i C _y 
\end{pmatrix} ,\\
\label{e4muminus} 
\hat\nue_- &
= -\frac{4  [x- x _k - i (y - y _k )]}{f |p-p _k |^4 }
 \begin{pmatrix}
 C _x - i C _y  \\
-C _\rho  
\end{pmatrix} ,\\
\label{e4mu3} 
\hat\nue _3 &=
 - \frac{2  }{f |p-p _k |^4 } 
  \left[[x- x _k - i (y - y _k ) ] 
  \begin{pmatrix}
C _\rho\\
C_ x + i C _y 
\end{pmatrix}  + \rho \begin{pmatrix}
C _x - i C _y \\
-C _\rho 
\end{pmatrix}  \right] .
\end{align} 

To obtain an eigenspinor $\nu \in E _\xi  ^- $ on $H ^3 $ it only remains to multiply $\hat \nue $ by the right power of the conformal factor, that is
\begin{equation}
\label{h3nu} 
\nu = \rho ^{ 3/2 } \hat \mu .
\end{equation}  

\section{Computation of $\omega (\nu , \Ispin (\thickbar\nu ) )  $  and $\omega (\nu , F \cdot \psi  ) $}
\label{integrals}
In this appendix we explain how the integrals  $\omega (\nu , \Ispin (\thickbar\nu ) )  $ and $\omega (\nu , F \cdot \psi  ) $ are evaluated explicitly.   Both integrals take the form
\begin{equation}
I = \int_{0}^\infty\int_{-\infty}^\infty\int_{\infty}^\infty\frac{\rho\, P_1(x^2,y^2,\rho^2,x_0^2)\,dx\,dy\,d\rho}{x_0^{4}[x_0^{2}x^4+(y^2+\rho^2+x_0^2)(x_0^2y^2+x_0^2\rho^2+1)+x^2(2x_0^2y^2+2x_0^2\rho^2+x_0^4-3)]^4},
\end{equation}
in which $P_1$ is a polynomial, so we will discuss them together. In the case of  $\omega (\nu ,  \Ispin (\thickbar\nu ) )  $ the polynomial $P _1 $ is given by
\begin{equation}
\begin{split}  
-16 (x_0^4-1)^2 
\Bigl[&
x^8 x_0^4
+2x^6x_0^2\left( 7-9x_0^4+2x_0^2(y^2+\rho^2) \right)\\ 
+&x^4 \Bigl(57+30x_0^2(y^2+\rho^2) +x_0^4\bigl(-44-7x_0^4-34x_0^2(y^2+\rho^2)+6(y^2+\rho^2)^2\bigr)\Bigr)\\ 
+2&x^2\Bigl(13(y^2+\rho^2)-3x_0^8(y^2+\rho^2)+2x_0^4(y^2+\rho^2)\bigl(-6+(y^2+\rho^2)^2\bigr)\\
 &\quad\quad   -x_0^6\bigl(9+7(y^2+\rho^2)^2\bigr)+x_0^2\bigl(7+9(y^2+\rho^2)^2\bigr)\Bigr)\\ 
 +&\left(x_0^2+y^2+\rho^2\right)^2 \left(1+x_0^2(y^2+\rho^2)\right)^2
\Bigr].
\end{split} 
\end{equation} 
In the case of  $\omega (\nu , F \cdot \psi  )  $ the polynomial $P _1 $ is given by
\begin{equation}
\begin{split}  
- 8 i  x _0 ^2 (x_0^4-1)^2  \Bigl[&
x^8x_0^2
+x^6\left(-30-18x_0^4+4x_0^2(y^2+\rho^2) \right)\\ 
-&x^4\Bigl(62(y^2+\rho^2)+x_0^2\bigl(95+7x_0^4+34x_0^2(y^2+\rho^2)-6(y^2+\rho^2)^2\bigr)\Bigr)\\ 
+2&x^2\Bigl(-3x_0^6(y^2+\rho^2)-17(y^2+\rho^2)^2+x_0^2(y^2+\rho^2)\bigl(-43+2(y^2+\rho^2)^2\bigr)\\ & \quad \quad  -x_0^4\bigl(20+7(y^2+\rho^2)^2\bigr)\Bigr)\\ &
+\left(x_0^2+y^2+\rho^2\right)^2\Bigl(-2(y^2+\rho^2)+x_0^2\bigl(-3+(y^2+\rho^2)^2\bigr)\Bigr)
\Bigr].
\end{split} 
\end{equation} 
The parameter $x_0$ satisfies $x_0\geq\sqrt{3}$.   

Since the integrand is even in both $x$ and $y$, we can restrict the ranges of integration to $0<x<\infty$ and $0<y<\infty$. We note that the denominator depends on $y,\rho$ only through the combination $\rho^2+y^2$.  So we substitute $(\rho,y)=(r\cos t,r\sin t)$ with $r>0$ and $0<t<\pi/2$.  Then the denominator is independent of $t$ and the numerator is a trigonometric polynomial in $t$.  So we can perform the $t$ integral, obtaining an expression of the form
\begin{equation}
I = \int_{0}^\infty\int_{0}^\infty\frac{P_2(x^2,r^2,x_0^2)dr\,dx}{x_0^4[x_0^2x^4+(r^2+x_0^2)(x_0^2r^2+1)+x^2(2x_0^2r^2+x_0^4-3)]^4},
\end{equation}
in which $P_2$ is a polynomial of degree 5 in $r^2$ and $4$ in $x^2$.

Next we try to integrate in $r$.  Since the integrand is an even function of $r$, we can extend the range of integration to $-\infty<r<\infty$.  Then we can evaluate the integral using complex analysis.  The degree of the numerator in $r$ is $10$ and the degree of the denominator in $r$ is $16$, so we can deform the region of integration away from the real axis to a closed contour in $\mathbb{C}$.  Then the integral can be evaluated by taking residues.

The poles are the zeros of the following quadratic in $r^2$:
\begin{equation}\label{quadratic}
r^4 + \left(\frac{x_0^4+1}{x_0^2}+2x^2\right)r^2+x^4+1+\frac{x_0^4-3}{x_0^2}x^2.
\end{equation}
This has solutions $\pm r_+$, $\pm r_-$, where $r_\pm^2$ are roots of the quadratic.  We find that $r_\pm^2$ are negative for all $x\geq0$ and $x_0\geq\sqrt{3}$, so $r_\pm$ can be chosen to lie on the positive imaginary axis.  Explicitly, they are given by
\begin{equation}
r_\pm = i\sqrt{x^2+\frac{x_0^4+1}{2x_0^2}\pm\sqrt{4\frac{x^2}{x_0^2}+\left(\frac{x_0^4-1}{2x_0^2}\right)^2}}.
\end{equation}

Having identified the poles, we can compute the integral in $r$.  We get two different expressions depending on whether we deform the contour to enclose the poles above or below the real axis:
\begin{align}
&\int_0^\infty \frac{P_2(x^2,r^2,x_0^2)}{x_0^4[x_0^2x^4+(r^2+x_0^2)(x_0^2r^2+1)+x^2(2x_0^2r^2+x_0^4-3)]^4}dr \nonumber\\
&=\frac{1}{2x_0^{12}}\left(
\Res_{r=r_+}\frac{P_2(x^2,r^2,x_0^2)}{(r-r_+^2)(r-r_-^2)}
+\Res_{r=r_-}\frac{P_2(x^2,r^2,x_0^2)}{(r-r_+^2)(r-r_-^2)}\right)\label{contour integral 1} \\
&=\frac{-1}{2x_0^{12}}\left(
\Res_{r=-r_+}\frac{P_2(x^2,r^2,x_0^2)}{(r-r_+^2)(r-r_-^2)}
+\Res_{r=-r_-}\frac{P_2(x^2,r^2,x_0^2)}{(r-r_+^2)(r-r_-^2)}\right)\label{contour integral 2}
\end{align}
The RHS of \eqref{contour integral 1} is a rational expression in $x^2,x_0^2,r_+,r_-$.  It is clearly invariant under permutation of $r_+,r_-$, so it can be rewritten in terms of $x^2$, $x_0^2$, $u=-i(r_++r_-)>0$ and $v=-r_+r_->0$.  Thus
\begin{equation}\label{rational integral}
I=\int_{0}^\infty R_1(x^2,x_0^2,u,v)dx
\end{equation}
where $R_1$ is the rational function defined by \eqref{contour integral 1}.  To determine how $u,v$ depend on $x$, we first note that the polynomial \eqref{quadratic} is equal to
\begin{equation}
(r+r_+)(r+r_-)(r-r_+)(r-r_-)=
(r^2+iur-v)(r^2-iur-v)=r^4+(u^2-2v)r^2+v^2.
\end{equation}
Equating coefficients with \eqref{quadratic} shows that
\begin{align}
v^2&=x^4+1+\frac{x_0^4-3}{x_0^2}x^2,\\
u^2-2v&=\frac{x_0^4+1}{x_0^2}+2x^2.
\end{align}
The first equation determines an elliptic curve.  The second determines a branched double cover over this elliptic curve, with the covering map being $(u,v,x)\mapsto (v,x)$.  There are two branch points and so the curve has genus 2.  This is too large to reduce the integral \eqref{rational integral} to elliptic integrals.

However, from the two expressions \eqref{contour integral 1} and \eqref{contour integral 2}, we see that $R_1$ changes sign when we replace $r_+,r_-$ with $-r_+,-r_-$.  So 
\begin{equation}\label{R1 sign change}
R_1(x^2,x_0^2,-u,v)=-R_1(x^2,x_0^2,u,v).
\end{equation}
This means that $R_1\,dx$ is invariant under the transformation $(x,u,v)\mapsto (-x,-u,v)$.  So we can rewrite $R_1\,dx$ in terms of the variables $(X,U,V)=(x^2,ux,v)$.  Doing so results in an expression
\begin{equation}
I=\int_{0}^\infty R_2(X,x_0^2,U,V)\,dX
\end{equation}
in which $R_2$ is rational.  The variables $X,U,V$ are related by
\begin{align}
V^2&=X^2+1+\frac{x_0^4-3}{x_0^2}X\label{curve1},\\
U^2&=X\left(\frac{x_0^4+1}{x_0^2}+2X+2V\right)\label{curve2}.
\end{align}
The first of these determines a conic, and the second is a branched double cover over this conic with four branch points.  It is therefore an elliptic curve.  So the integral $I$ can be reduced to elliptic integrals.

To evaluate the integral, we first parametrise the conic \eqref{curve1} with a single variable $t$.  A convenient rational parametrisation is
\begin{align}
X&=\frac{1-t}{4tx_0^2}\big((x_0^2+3)(x_0^2-1)-(x_0^2-3)(x_0^2+1)t\big)\\
V&=\frac{(x_0^2+3)(x_0^2-1)}{4x_0^2t}-\frac{(x_0^2-3)(x_0^2+1)}{4x_0^2}t.
\end{align}
Then the integration limits $(X,V)=(0,1)$ and $(X,V)\to(+\infty,+\infty)$ correspond to $t=1$ and $t\to0^+$ respectively.  Moreover, it can be checked that the integration measure becomes $dX= -V dt/t$.  Inserting this into \eqref{curve2} gives
\begin{equation}\label{elliptic curve}
w^2=(1-t)\big((x_0^2+3)(x_0^2-1)-(x_0^2-3)(x_0^2+1)t\big)\big((x_0^2+3)(x_0^2-1)+4t\big)
\end{equation}
in which $w=2x_0^2Ut$.  As has already been observed, this is an elliptic curve.  The integral $I$ reduces to an elliptic integral of the form
\begin{equation}
I=\int_0^1 R_3(w^2,t) \frac{dt}{w}
\end{equation}
for some rational function $R_3$.  We have written $R_3$ as a function of $w^2$, rather than $w$, because from \eqref{R1 sign change} the integrand $R_3\, dt/w$ changes sign under the transformation $(w,t)\mapsto(-w,t)$ and hence $R_3$ is invariant under this transformation.  On substituting for $w^2$ using \eqref{elliptic curve}, $R_3$ becomes a rational function of $t$.  The numerator and denominator of this rational function are polynomials in $t$ of degrees 6 and 11 respectively.

To reduce the integral to a standard form, we eliminate the poles from $R_3(t)$.  We can find a rational function $Q(t)$ and a degree one polynomial $L(t)$ such that
\begin{equation}
R_3\frac{dt}{w}=d(wQ)+L\frac{dt}{w}.
\end{equation}
Then, since $w (1) =0$,
\begin{equation}
I=-Q(0)w(0)+\int_0^1 L(t)\frac{dt}{w}.
\end{equation}
Finally, to evaluate the integral of $L\,dt/w$, we make the substitution
\begin{equation}
t = \frac{(x_0^2+1)^2}{4}(T^2-1)+1,\quad w=\frac{(x_0^2+1)^2}{4}WT.
\end{equation}
Then
\begin{equation}
dt = \frac{(x_0^2+1)^2}{2}T\,dT,\quad W^2=(1-T^2)((x_0^2+3)(x_0^2-1)^3-(x_0^2-3)(x_0^2+1)^3T^2).
\end{equation}
The integral becomes
\begin{equation}
I=-Q(0)w(0)+\int_{\sin\varphi}^1 \frac{L\left((x_0^2+1)^2(T^2-1)/4+1\right)}{\sqrt{(x_0^2+3)(x_0^2-1)^3}}\frac{dT}{\sqrt{(1-T^2)(1-mT^2)}},
\end{equation}
in which
\begin{equation}
\sin^2\varphi = 1-\frac{4}{(x_0^2+1)^2},\quad m=\frac{(x_0^2-3)(x_0^2+1)^3}{(x_0^2+3)(x_0^2-1)^3}.
\end{equation}
Since $L$ has degree 1, this integral is a linear combination of elliptic integrals of the first and second kinds with parameter $m$.  The final expressions \eqref{intnuibnufinal}, \eqref{intnunukillfinal} that we obtain for the integrals are remarkably simple, given the complexity of the functions being integrated.

Finally, it is interesting to compare the elliptic curve \eqref{elliptic curve} that underlies this calculation with the monopole spectral curve.  Completing the square brings the spectral curve \eqref{2-monopole curve} to the form
\begin{equation}
\label{2-monopole curve square} 
W^2=z^4+\frac{x_0^8-6x_0^4-3}{4x_0^2}z^2+1,
\end{equation}
in which $W=(x_0^2-z^2)w/x_0+z(1-x_0^4)/(2x_0)$.  This is a branched double cover over a complex line with branch points at the roots $\pm z_+,\pm z_-$ of the quartic (\ref{2-monopole curve square}).  These roots satisfy
\begin{equation}
z_+^2+z_-^2=-\frac{x_0^8-6x_0^4-3}{4x_0^2},\quad z_+^2z_-^2=1,
\end{equation}
and we can choose to label them so that $z_+z_-=1$.  Their cross ratio is therefore
\begin{equation}
(-z_-,z_-;-z_+,z_+)=\frac{(z_+-z_-)^2}{(z_++z_-)^2}=\frac{z_+^2+z_-^2-2z_+z_-}{z_+^2+z_-^2+2z_+z_-}=\frac{(x_0^2+3)(x_0^2-1)^3}{(x_0^2-3)(x_0^2+1)^3}.
\end{equation}
On the other hand, the cross ratio of the branch points of \eqref{elliptic curve} is
\begin{equation}
\left(\infty,-\frac{(x_0^2+3)(x_0^2-1)}{4};1,\frac{(x_0^2+3)(x_0^2-1)}{(x_0^2-3)(x_0^2+1)} \right)=\frac{(x_0^2+3)(x_0^2-1)^3}{(x_0^2-3)(x_0^2+1)^3}.
\end{equation}
Since the two cross-ratios of the branch points are equal, the $j$-invariants are equal and the curve \eqref{elliptic curve} that appears in our elliptic integrals is biholomorphic to the monopole spectral curve \eqref{2-monopole curve}.

\bibliographystyle{amsplain}
\bibliography{biblio}

\end{document}